# TRULY HYPERCOMPLEX NUMBERS:

# UNIFICATION OF NUMBERS AND VECTORS

Redouane BOUHENNACHE
(pronounce Redwan Boohennash)
Independent Exploration Geophysical Engineer / Geophysicist
14, rue du 1$^{er}$ Novembre,
Beni-Guecha Centre, 43019 Wilaya de Mila,
Algeria
E-mail: redouane.bouhennache@outlook.com

First written: 21 July 2014
Revised: 17 May 2015

## Abstract

Since the beginning of the quest of hypercomplex numbers in the late eighteenth century, many hypercomplex number systems have been proposed but none of them succeeded in extending the concept of complex numbers to higher dimensions. This paper provides a definitive solution to this problem by defining the truly hypercomplex numbers of dimension $N \geq 3$. The secret lies in the definition of the multiplicative law and its properties. This law is based on spherical and hyperspherical coordinates. These numbers which I call spherical and hyperspherical hypercomplex numbers define Abelian groups over addition and multiplication. Nevertheless, the multiplicative law generally does not distribute over addition, thus the set of these numbers equipped with addition and multiplication does not form a mathematical field. However, such numbers are expected to have a tremendous utility in mathematics and in science in general.

## Note

This paper (or say preprint or e-print) has been submitted, under the title "Spherical and Hyperspherical Hypercomplex Numbers: Merging Numbers and Vectors into Just One Mathematical Entity", to the following journals:
- *Bulletin of Mathematical Sciences* on 08 August 2014,
- *Hypercomplex Numbers in Geometry and Physics* (HNGP) on 13 August 2014 and has been accepted for publication on 29 April 2015 in issue No. 22 of HNGP.

# 1. Introduction

I discovered the spherical and hyperspherical hypercomplex numbers before dawn on a night of the holy Muslim fasting month Ramadan of the year 1424 of Hegira, corresponding to mid-November 2003. I considered the problem of extending the complex numbers to higher dimensions starting from dimension 3, and in a quarter an hour I figured out the solution. An abstract of my article has been published in the proceedings of "Number, time and relativity conference" in Moscow in 2004 under the name spelling Bouhannache[1]. First I thought those numbers formed a field, because I did not conduct a thorough check of all the field axioms. Later on 25 March 2006, I realized that the multiplicative law was not distributive on addition, making these numbers form Abelian groups only. I did not have enough time to finish my article because I was working in the Algerian desert from 2004 to 2005. And from 2006 till 2013 I was busy. I resumed writing my article in November 2013.

# 2. Brief history of the hypercomplex numbers' quest

Quaternions were discovered by William Rowan Hamilton in 1843. Hamilton was looking for ways of extending complex numbers (which can be viewed as points on a plane) to higher spatial dimensions. He could not do so for 3 dimensions, but 4 dimensions produce quaternions. According to a story he told, he was out walking one day with his wife when the solution in the form of equation $i^2 = j^2 = k^2 = ijk = -1$ suddenly occurred to him; he then promptly carved this equation into the side of nearby Brougham bridge (now called Broom Bridge) in Dublin.[2]

This involved abandoning the commutative law, a radical step for the time. Vector algebra and matrices were still in the future. Not only this, but Hamilton had in a sense invented the cross and dot products of vector algebra. Hamilton also described a quaternion as an ordered four-element multiple of real numbers, and described the first element as the 'scalar' part, and the remaining three as the 'vector' part. If two quaternions with zero scalar parts are multiplied, the scalar part of the product is the negative of the dot product of the vector parts, while the vector part of the product is the cross product. But the significance of these was still to be discovered.[3]

Quaternions have been discovered by Hamilton after 10 years of research of possible fields of triples or quadruples. They form a skew field noted H where multiplication is associative but not commutative.

Octonions are a non-associative extension of the quaternions. They were discovered by John T. Graves in 1843, and independently by Arthur Cayley, who published the first paper on them in 1845. They are sometimes referred to as Cayley numbers or the

---

[1] Redouane Bouhannache, Solving the old problem of hypercomplex number fields: The new commutative and associative hypercomplex number algebra and the new vector field algebra, Number, Time, Relativity. Proceedings of International Scientific Meeting, Moscow: 10 – 13 August, 2004 https://sites.google.com/site/bouhannache/home/files-containing-my-contributions-on-hypercomplex-numbers Accessed 05 June 2013, http://hypercomplex.xpsweb.com/articles/194/en/pdf/sbornik.pdf Accessed 05 June 2013
[2] Quaternion, http://www.fact-index.com/q/qu/quaternion.html Accessed 10 Jan 2004
[3] Ibid.

Cayley algebra.[4] They are 8-dimensional. Their multiplicative law is neither commutative nor associative but alternative.

Sedenions form a 16-dimensional algebra over the reals where multiplication is neither commutative, nor associative, nor alternative but power associative.

Quaternions, Octonions, sedenions, etc. are defined only for dimensions $N=2^m$, $m \geq 2$ and are obtained through the Cayley-Dickson construction.

About mid 19th century, Grassmann lies the foundations of exterior algebra. In 1879 Clifford publishes "Applications of Grassmann's extensive algebra" formulating what is now known as Clifford algebra.[5] The latter is a geometric algebra. He introduced a new type of hybrid product called « the geometric product », which is associative, non-commutative and of a heterogeneous nature.

A lot of other hypercomplex number systems have been proposed but none of them form an extension of the complex numbers to higher dimensions. However, the hypercomplex number system I discovered does.

In a letter Hamilton wrote to his son Archibald, he says:

« Every morning in the early part of the above-cited month [October 1843] on my coming down to breakfast, your brother William Edwin and yourself used to ask me, 'Well, Papa, can you multiply triples [triplets]?' Whereto I was always obliged to reply, with a sad shake of the head, 'No, I can only add and subtract them.' »[6, 7]

This paper gives the answer to the question Hamilton's sons asked their father, showing not only how to multiply and divide triples but also how to do it for any *n*-tuples (tuplets).

I call the hyperspherical hypercomplex numbers the true hypercomplex numbers or just hypercomplex numbers.

## 2.1. Epistemological reasons standing behind mathematicians' failure to find the hyperspherical hypercomplex numbers - Main error of the previous attempts

Most mathematicians' researches (including Hamilton, Dickson, Cayley) have been misled by considering that, like in the case of complex numbers, there must exist a linear relation between the coordinates (components) of the hypercomplex product $(x_p, y_p, z_p,...)$ and the individual number coordinates $(x, y, z,...)$ and $(x', y', z',...)$. Such a relation doesn't exist at all.

---

[4] Octonion, Wikipedia, the free encyclopedia http://en.wikipedia.org/wiki/Octonion Accessed 03 April 2004

[5] Sketching the History of Hypercomplex Numbers, http://history.hyperjeff.net/hypercomplex Accessed 03 April 2004

[6] History of quaternions, Wikipedia, the free encyclopedia, http://en.wikipedia.org/wiki/History_of_quaternions Accessed 13 March 2014

[7] A Brief History of Quaternions, http://world.std.com/~sweetser/quaternions/intro/history/history.html Accessed 03 April 2004

Hence most of them have been using linear algebra (matrices) to solve that problem, whereas the latter can never be solved using conventional matrices (expressing linear applications).

# 3. Construction of the spherical hypercomplex numbers: The new concept

## 3.1. Introduction: Position in 3D space – Spherical coordinates

Spherical coordinates have coordinates typically named $r$, $\theta$ and $\varphi$ where the radius $r$ ranges from 0 to $\infty$, the azimuth $\theta$ ranges from 0 to $2\pi$, and the colatitude $\phi$ ranges from 0 to $\pi$. They describe a point in space as follows: from the origin (0,0,0), go $r$ units along the $z$-axis, rotate $\phi$ down from the $z$-axis in the $x$-$z$ plane (colatitude), and rotate $\theta$ counterclockwise about the $z$-axis (azimuth or longitude).[8] The name of the system comes from the fact that the simple equation $r = 1$ describes the unit sphere.

Unlike the usual system which uses the co-latitude, I will use the latitude $\varphi = \dfrac{\pi}{2} - \phi$ like in geographic coordinate system.

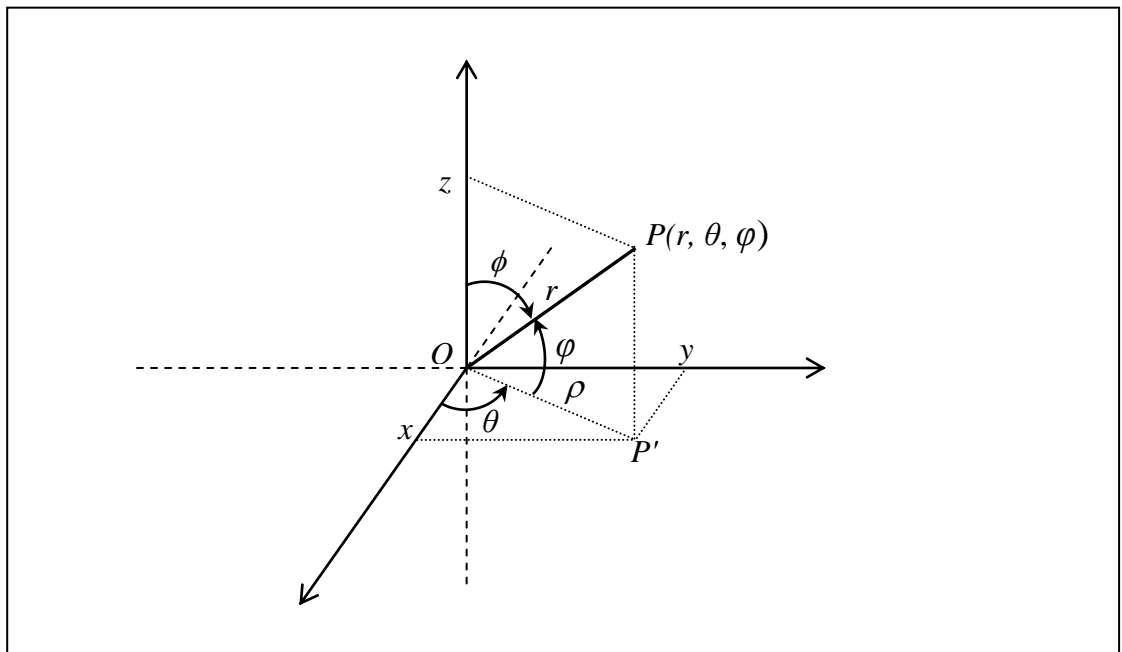

Fig. 3.1. *Spherical coordinates*

$\varphi$ ranges from $-\dfrac{\pi}{2}$ to $+\dfrac{\pi}{2}$

A point $P(x,y,z)$ coordinates (or say vector $\overrightarrow{OP}$ components) are:

$$
\begin{aligned}
x &= r\cos\varphi\cos\theta \\
y &= r\cos\varphi\sin\theta \\
z &= r\sin\varphi
\end{aligned}
\qquad (3.1)
$$

---

[8] Spherical coordinate system, Wikipedia, the free encyclopedia
http://en.wikipedia.org/wiki/Spherical_coordinate_system  Accessed 25 November 2003

## 3.2. Key principle of the spherical and hyperspherical hypercomplex numbers' construction: Orthogonal-complex-plane construction

I will introduce hereafter a new concept in constructing hypercomplex numbers. It employs two orthogonal complex planes. I therefore call it the orthogonal-complex-plane construction.

### 3.2.1. Original idea standing behind my discovery of hypercomplex construction (Nov 2003). Physical and geometrical meaning. Cartesian construction of 3D hypercomplex numbers

In fact, I originally considered the following:

In spherical coordinates, in addition to the angle $\theta$, we have a new angle $\varphi$. Let's consider the plane containing the point $P(r,\theta,\varphi)$ and $z$-axis (hence containing also the origin $O$) and which is perpendicular (or normal) to the plane $(xOy)$. Let's limit our imagination to that plane only – i.e. not taking into account the angle $\theta$. In such a plane, $P$'s Cartesian coordinates will be $\rho$ (abscissa) and $z$ (ordinate). $P'$ is the orthogonal projection of $P$ on the plane $(xOy)$. See Fig.3.1.

$\rho = r\cos\varphi$ is the 2D conventional complex modulus.

It is obvious that the position of any point $P$ can also be determined in this plane $(P'Oz)$ with polar coordinates $r$ and $\varphi$. But we should note here that, in this plane, $\varphi$ will have any value, and not necessarily comprised between $-\frac{\pi}{2}$ and $+\frac{\pi}{2}$ as in the definition of the spherical latitude.

Similarly to the complex plane $(xOy)$, we certainly can consider the new plane $(P'Oz)$ as being complex. The argument here will be $\varphi$ and the modulus $r$.

Let's consider we have two 3D hypercomplex numbers $h$ and $h'$
$$h = x + iy + jz, \; h' = x' + iy' + jz' \qquad (3.2)$$

Assume there exists a certain 3-D commutative and associative multiplication law. This law must verify the following 2 conditions:
1. It must preserve 2-D $(xOy)$ complex plane arithmetics when $z=z'=0$.
2. Complex arithmetic laws must hold on a new complex plane $\pi(\rho,z)$ which is insensitive to 2-dimentional $\theta$ (insensitive to rotation on $\pi(x,y)$ plane, that is to say $(xOy)$ plane) and has only argument $\varphi$. The complex plane $(P'Oz)$ contains the points $P$ and $P'$ and is normal to the conventional complex plane $(xOy)$.
The abscissas will be $\rho$, the module of $(x,y)$ component vectors (exactly equivalent to modulus of complex number $c=x+iy$). $\rho = \sqrt{x^2 + y^2}$, $\rho' = \sqrt{x'^2 + y'^2}$

Now let's limit our imagination to the $(P'Oz)$ plane which is invariant to $\theta$ and $\rho$ is always positive by definition ($\rho \geq 0$). We will have a 2-D complex plane with imaginary number $j$ (the same laws of conventional complex plane equipped with imaginary number *i* will hold for the new plane equipped with *j*). The two complex numbers in question can be written in this new plane as follows:
$h = \rho + jz$ and $h' = \rho' + jz'$ (3.3)

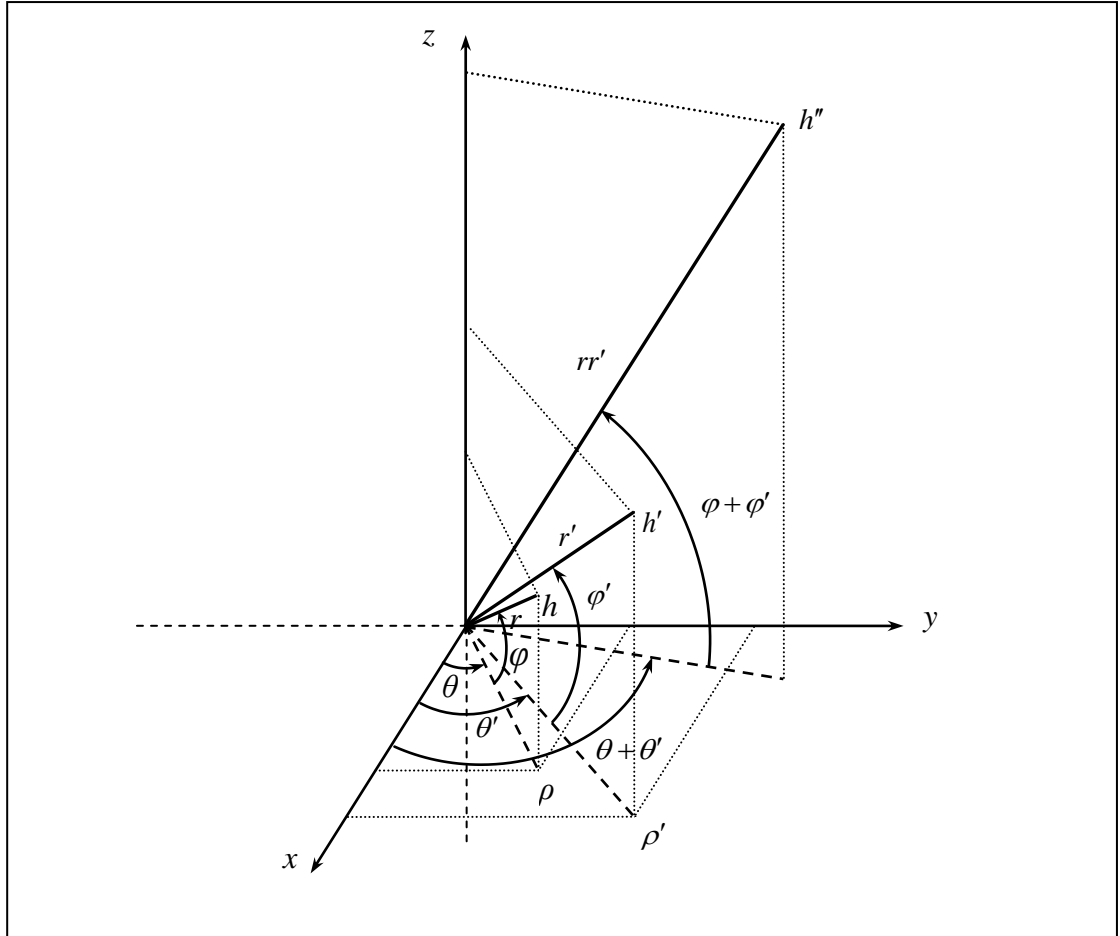

Fig.3.2. *Multiplication of two triples h and h'.*

When multiplying them, arguments $\varphi$'s will add up. In Cartesian coordinates we get:
Product: $h'' = (\rho\rho' - zz') + j(\rho z' + \rho' z)$ (3.4)

This means we get a *(P'Oz)*-plane-complex number whose modulus is *rr'* and whose real part is $\rho\rho' - zz'$. And we see that the *z*-component is $\rho z' + \rho' z$. See Fig.3.3.

After conducting the multiplication in *(P'Oz)* plane, we multiply the abscissa and ordinate by $\cos(\theta + \theta')$ and $\sin(\theta + \theta')$ respectively, because arguments sum up in the conventional complex plane *(xOy)*.

$$\cos(\theta + \theta') = \cos\theta\cos\theta' - \sin\theta\sin\theta' = \frac{x}{\rho}\frac{x'}{\rho'} - \frac{y}{\rho}\frac{y'}{\rho'} = \frac{xx' - yy'}{\rho\rho'}$$ (3.5)

$$\sin(\theta + \theta') = \sin\theta\cos\theta' + \cos\theta\sin\theta' = \frac{y}{\rho}\frac{x'}{\rho'} + \frac{x}{\rho}\frac{y'}{\rho'} = \frac{xy' + x'y}{\rho\rho'}$$ (3.6)

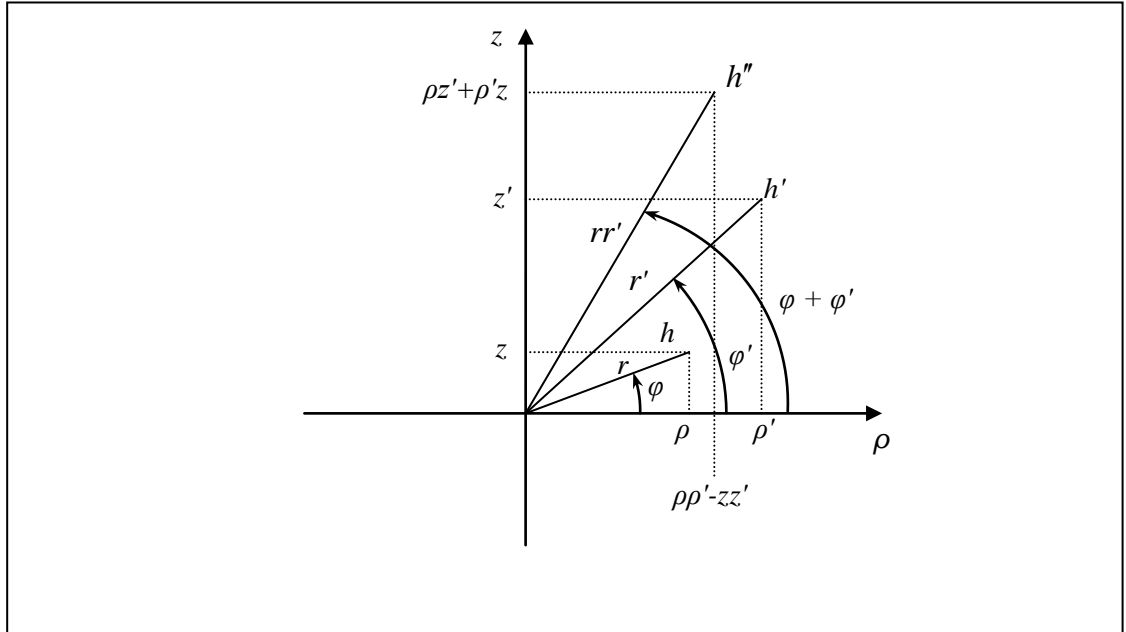

Fig.3.3. *Multiplication in the orthogonal planes containing Oz*

Hence

$$x_{(product)} = x'' = \frac{(xx' - yy')}{\rho\rho'}(\rho\rho' - zz') \tag{3.7}$$

$$y_{(product)} = y'' = \frac{(xy' + x'y)}{\rho\rho'}(\rho\rho' - zz') \tag{3.8}$$

$$x_{(product)} = x'' = (xx' - yy')\left(1 - \frac{zz'}{\rho\rho'}\right) \tag{3.9}$$

$$y_{(product)} = y'' = (xy' + x'y)\left(1 - \frac{zz'}{\rho\rho'}\right) \tag{3.10}$$

Note: $1 - \dfrac{zz'}{\rho\rho'} = 1 - \tan\varphi\tan\varphi'$ (3.11)

Hence, as a recap, $h = (x, y, z)$ and $h' = (x', y', z')$ are two 3D hypercomplex numbers. $h'' = h \cdot h' = (x'', y'', z'')$ is their product. The spherical hypercomplex multiplication is defined in Cartesian coordinates by:

$$\begin{cases} x'' = (xx' - yy')\left(1 - \dfrac{zz'}{\sqrt{x^2 + y^2}\sqrt{x'^2 + y'^2}}\right) \\ y'' = (xy' + x'y)\left(1 - \dfrac{zz'}{\sqrt{x^2 + y^2}\sqrt{x'^2 + y'^2}}\right) \\ z'' = z\sqrt{x'^2 + y'^2} + z'\sqrt{x^2 + y^2} \end{cases} \tag{3.12}$$

$$\rho = \sqrt{x^2 + y^2},\ \rho' = \sqrt{x'^2 + y'^2}. \tag{3.13}$$

This is valid for $\rho$ and $\rho'$ different from zero. If one or the two of them equal zero, then x, y, and z can be calculated using the geometric form (See paragraph 4.1.).

If $\rho = 0$ and $\rho' = 0$, then $\rho\rho' - zz' = -zz'$ and $\theta$ and $\theta'$ must be provided for $h$ and $h'$, that is to say the longitudes of planes containing $h$ and $h'$ and which are normal to complex plane must be provided.

$x_{(product)} = x'' = -zz'\cos(\theta + \theta')$ (3.14)

$y_{(product)} = y'' = -zz'\sin(\theta + \theta')$ (3.15)

If $\rho = 0$ and $\rho' \neq 0$ then $\rho\rho' - zz' = -zz'$ and $\theta$ must be provided for $h$:

$$x'' = -zz'(\cos\theta\cos\theta' - \sin\theta\sin\theta') = -zz'(\cos\theta\frac{x'}{\rho'} - \sin\theta\frac{y'}{\rho'})$$

$$x'' = \frac{zz'}{\rho'}(y'\sin\theta - x'\cos\theta)$$

(3.16)

$$y'' = -zz'(\sin\theta\cos\theta' + \cos\theta\sin\theta') = -zz'(\sin\theta\frac{x'}{\rho'} + \cos\theta\frac{y'}{\rho'})$$

$$y = -\frac{zz'}{\rho'}(x'\sin\theta + y'\cos\theta)$$

(3.16)

The name spherical hypercomplex number applies only to 3D hypercomplex numbers and derives from the fact that these numbers and their multiplicative operation are defined using spherical coordinates.

The name hyperspherical hypercomplex numbers applies to numbers with dimension $N \geq 4$ and derives from the fact that these numbers and their multiplicative operation are defined using hyperspherical coordinates.

Or if we consider the *n*-hypersphere to be an *n*-sphere and the hyperspherical coordinates are spherical then we call all these numbers simply "spherical hypercomplex numbers".

### 3.2.2. Second method

We consider the assumptions made in section 3.2.1.

Similarly to the complex plane *(xOy)*, we certainly can consider the new plane *(P'Oz)* as being complex. The argument here will be $\varphi$ and the modulus *r*. If P is located on the peripheral of a unit circle (radius = number modulus= $r = 1$), the famous classical Euler's formula will also be valid in the new orthogonal plane *(P'Oz)*.

Hence it satisfies:

$e^{j\varphi} = \cos\varphi + j\sin\varphi$    w h e r e    $j^2 = -1$    or    $j = +\sqrt{-1} = j = e^{j\frac{\pi}{2}}$    and
$j^2 = e^{j\pi} = \cos\pi = -1$

Remember that this is valid only in the plane *(P'Oz)* normal to *(xOy)*.

*j* is the unit imaginary number of *z*-axis in the new complex plane *(P'Oz)* (equivalent: unit vector of *z*-axis).

Therefore any complex number $c'$ of modulus $r$ and argument $\varphi$ can be expressed in the new plane, in the conventional way, as follows:

$$c' = re^{j\varphi} = r(\cos\varphi + j\sin\varphi) = r\cos\varphi + jr\sin\varphi \tag{3.17}$$

The conventional *(x,y)* complex number component can be denoted or regarded as being $c = x + iy = r\cos\varphi\cos\theta + ir\cos\varphi\sin\theta$ (3.18)

### 3.2.3. How to link the two orthogonal complex planes *(xOy)* and *(P'Oz)*

### 3.2.3.a. First Approach: New 3-D Euler's Formula

Is there a way to merge the two notations (3.17) and (3.18) into just one? Can we explicitly express point *P* coordinates *(x,y,z)* in function of the two Euler's formulas so that we get the conventional expression of *(x,y,z)* in function of spherical coordinates *(r,θ,φ)*?

Let's consider $h$ the 3-D hypercomplex number expressing the position of *P*.
In the plane *(P'Oz)*, $c'$ expression will be:

$$c' = re^{j\varphi} = r\cos\varphi + jr\sin\varphi = \rho + jr\sin\varphi \tag{3.19}$$

On the conventional complex plane *(xOy)* we have:
$$c = x + iy = \rho e^{i\theta} \tag{3.20}$$

By definition of *(x,y,z)* Cartesian coordinates in function of spherical ones, $h$ coordinates are:
$$\begin{aligned} x &= r\cos\varphi\cos\theta \\ y &= r\cos\varphi\sin\theta \\ z &= r\sin\varphi \end{aligned} \tag{3.21}$$

$\Rightarrow h = x + iy + jz = r\cos\varphi\cos\theta + ir\cos\varphi\sin\theta + jr\sin\varphi = r\cos\varphi(\cos\theta + i\sin\theta) + jr\sin\varphi$

$$h = r\cos\varphi\, e^{i\theta} + jr\sin\varphi \tag{3.22}$$

Let's go back to definition of the plane *(P'Oz)*. The unit of abscissas $\rho$ equals *1* in this plane:
$$unit(\rho) \cdot j = 1 \cdot j = j \cdot 1 = j \tag{3.23}$$

What does this unit represent in plane *(xOy)*? It represents a complex number of modulus $\sqrt{x^2 + y^2} = 1$ and argument $\theta$. We can link the two formulas by assuming that:
$$e^{i\theta} j = je^{i\theta} = j \tag{3.24}$$

Replacing $j$ in equation (3.22) we can write:
$h = r\cos\varphi\, e^{i\theta} + je^{i\theta} r\sin\varphi = re^{i\theta}(\cos\varphi + k\sin\varphi)$
Hence:
$$h = re^{i\theta} e^{j\varphi} = re^{i\theta + j\varphi} \tag{3.25}$$

This is the analytical expression of P position in 3-D space in function of its spherical coordinates $r$, $\theta$ and $\varphi$.

Thus we can express $h$ in function of spherical coordinates using Euler's formulas for $\theta$ and $\varphi$.

Equation (3.25) is called the *geometric form* (or *exponential form*) of the 3D hypercomplex number $h$.

For the particular case of $r=1$ we get:
$$h = e^{i\theta} e^{j\varphi} = \cos\theta \cos\varphi + i\sin\theta \cos\varphi + j\sin\varphi \qquad (3.26)$$
This is what I call the 3-D or spherical Euler's formula. The classical 2-D formula is circular or polar.

$h = r\cos\varphi\cos\theta + ir\cos\varphi\sin\theta + jr\sin\varphi$ is called the *trigonometric form* of the 3D hypercomplex number $(r,\theta,\varphi)$.

### 3.2.3.b. Second approach

We can also adopt the inverse method:
In the *(P'Oz)*-plane

c' expression will be: $c' = re^{j\varphi}$ \qquad (3.27)

Let's assume that $h = re^{i\theta + j\varphi} = re^{i\theta}e^{j\varphi}$, developing, comparing (or identifying) with spherical-to-Cartesian transformation:
$$h = re^{i\theta}e^{j\varphi} = re^{i\theta}(\cos\varphi + j\sin\varphi) = re^{i\theta}\cos\varphi + re^{i\theta}j\sin\varphi \qquad (3.28)$$
$$h = r\cos\theta\cos\varphi + ir\sin\theta\cos\varphi + je^{i\theta}r\sin\varphi = x + iy + jz \qquad (3.29)$$
By definition $h$'s Cartesian coordinates are:
$$\begin{aligned} x &= r\cos\varphi\cos\theta \\ y &= r\cos\varphi\sin\theta \\ z &= r\sin\varphi \end{aligned} \qquad (3.30)$$
By comparison we deduce:
$$e^{i\theta} j = je^{i\theta} = j(\cos\theta + i\sin\theta) = j \qquad (3.31)$$

This means $j$ multiplied by any *(x,y)* complex number of modulus 1 (vector of module 1) remains unchanged whatever the argument $\theta$ is. Hence, $j$ is insensitive to rotation in the *(xOy)* classical complex plane.

> ***Theorem***: $j$ multiplied by a complex number of modulus 1 remains unchanged. $j$ is insensitive to rotation in the complex plane (or $\theta$ (2D) argument change.)

## 4. Geometric form of the spherical hypercomplex numbers

$h=(x, y, z)$ is a spherical hypercomplex number.
$|h| = r = \sqrt{x^2 + y^2 + z^2}$ is called the modulus of this number
$$\begin{aligned} x &= r\cos\varphi\cos\theta \\ y &= r\cos\varphi\sin\theta \\ z &= r\sin\varphi \end{aligned} \qquad (4.1)$$

$$\cos\theta = \frac{x}{\sqrt{x^2+y^2}}, \ \sin\theta = \frac{y}{\sqrt{x^2+y^2}} \tag{4.2}$$

Any $\theta$ satisfying this is called a *longitudinal argument* or simply *longitude* and denoted $\theta(h)$ or $\arg\theta(h)$.

The angle $\theta \in [0,+2\pi[$ satisfying this is called the *principal longitudinal argument* and denoted $\mathrm{Arg}\theta(h)$.

$$\arg\theta(h) = \mathrm{Arg}\theta(h) + 2k\pi, \ k \in \mathbb{Z} \tag{4.3}$$

$$\sin\varphi = \frac{z}{r} \tag{4.4}$$

$$\tan\varphi = \frac{z}{\rho} \tag{4.5}$$

$\rho = \sqrt{x^2+y^2}$ is the complex (x,y) modulus

Any $\varphi$ satisfying this and comprised in $\left[-\frac{\pi}{2},+\frac{\pi}{2}\right]+2k\pi$ is called a *first latitudinal argument* or *first latitude* denoted $\varphi(h)$ or $\arg\varphi(h)$

Any hypercomplex number, $h = x + yi + jz$ has infinitely many arguments but they all differ by multiples of $2\pi$.

*(r, θ,φ)* is the *geometric form* of the 3D-hypercomplex number *h*.
*(r, θ+π, π–φ)* is the *duplicate* (or *replicate*) *geometric form*. See Fig.4.1.
*(r, θ,φ)* and *(r, θ+π, π–φ)* have the same Cartesian coordinates.
We say that *(r, θ,φ)* and *(r, θ+π, π–φ)* are *duplicates* or *replicates*.

Due to spherical duplicity (duplication / replication) of the geometric form of the spherical hypercomplex number:

Any $\varphi$ satisfying equation (4.4) and comprised in $\left]+\frac{\pi}{2},+\frac{3\pi}{2}\right[+2k\pi$ is called a *second latitudinal argument* or *second latitude* and denoted $\varphi_2(h)$ or $\arg\varphi_2(h)$.

Also the principal first and second latitudes can be defined as:
$$\arg\varphi(h) = \mathrm{Arg}\varphi(h) + 2k\pi, k \in \mathbb{Z} \tag{4.6}$$
$\mathrm{Arg}\varphi_2(h)$ is always used along with $\mathrm{Arg}\theta(h)+\pi$. See Fig.4.1.
$$\mathrm{Arg}\varphi_2(h) = \pi - \mathrm{Arg}\varphi(h) \tag{4.7}$$

In this paper we consider or we note: $\theta = \mathrm{Arg}\theta(h)$ and $\varphi = \mathrm{Arg}\varphi(h)$ \hfill (4.8)

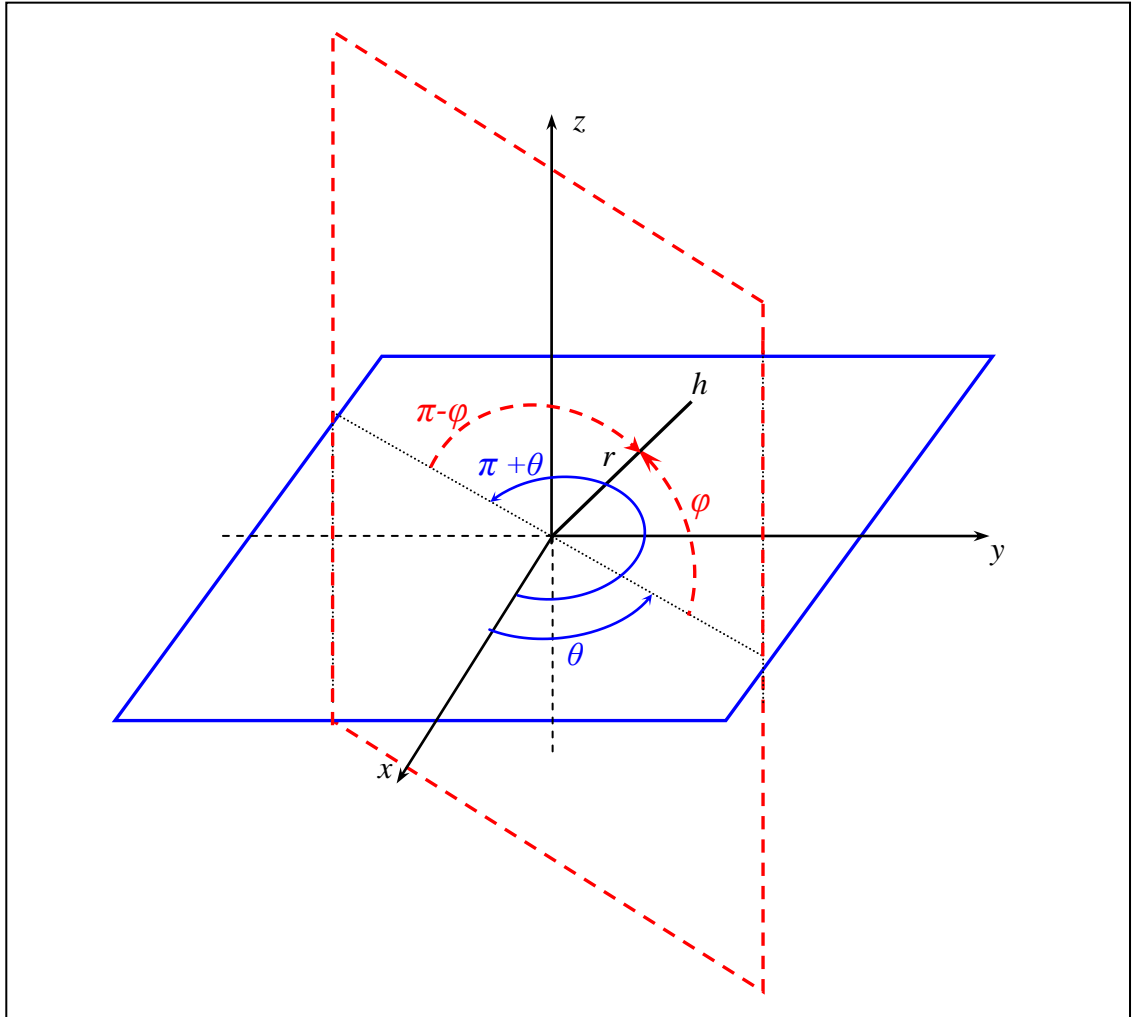

Fig.4.1. *Duplicity or duplication/replication of 3D spherical hypercomplex numbers*

***Important Note:*** Since hypercomplex numbers express a point position coordinates in spherical coordinates, or vector components in the same coordinate system, we can express such numbers this way: *h(r,θ,φ)* or just *(r,θ,φ)* or *r(θ,φ)*.

***Note:*** A hypercomplex number *h* can be expressed *(x,y,z)* or *h(x,y,z)* or using the geometric form *(r, θ, φ)* or *h(r, θ, φ)*.

## 4.1. Multiplication of the two spherical hypercomplex numbers *h* and *h'* using the geometric form

$h = r\, e^{i\theta + j\varphi}$      (4.9)
$h' = r'\, e^{i\theta' + j\varphi'}$      (4.10)
$hh' = rr'\, e^{i(\theta+\theta') + j(\varphi+\varphi')}$      (4.11)

When two spherical hypercomplex numbers are multiplied, the moduli are multiplied and the arguments added two by two.
if $(x, y) = (0,0)$ and/or $(x', y') = (0,0)$, then one of/the two hypercomplex numbers has/have to be determined by the angles $\theta$ and $\theta'$. It's not like in complex case where (0,0) has no determined $\theta$ and any $\theta$ satisfies the relationship.

When $\rho = 0$ that is to say $\varphi = \pm\frac{\pi}{2}$ the Cartesian coordinates formula (3.12) won't work. We use the geometric formula for the calculation of *x* and *y* components.

## 5. Is the set $\mathbb{H}$ of 3D hypercomplex numbers equipped with addition and multiplication a field?

For simplification reasons, we will use the geometric hyperspherical form of hypercomplex numbers. The reader can verify the field axioms using Cartesian expressions.

A field is a set *F* that is a commutative group with respect to two compatible operations, addition and multiplication, with "compatible" being formalized by distributivity, and the caveat that the additive identity (0) has no multiplicative inverse (one cannot divide by 0).

The most common way to formalize this is by defining a field as a set together with two operations, usually called addition and multiplication, and denoted by + and ·, respectively, such that the field axioms hold; subtraction and division are defined implicitly in terms of the inverse operations of addition and multiplication.[9]

Let's consider the set $\mathbb{H}$ of triples (*x*,*y*,*z*) such that *x*,*y* and *z* are real. This set is equipped with the two binary operations + and · defined for *h* and *h'* as follows:

$h + h' = (x, y, z) + (x', y', z') = (x + x', y + y', z + z')$ (5.1)

$h \cdot h' = \left( (xx' - yy')\left(1 - \frac{zz'}{\rho\rho'}\right), (xy' + x'y)\left(1 - \frac{zz'}{\rho\rho'}\right), \rho z' + \rho' z \right), \rho \neq 0, \rho' \neq 0$ (5.2)

Such that: $\rho = \sqrt{x^2 + y^2}$, $\rho' = \sqrt{x'^2 + y'^2}$ (5.3)

$h \cdot h'$ can be expressed in function of $(r, \theta, \varphi)$ and $(r', \theta', \varphi')$ such that:

$r = \sqrt{x^2 + y^2 + z^2}$, (5.4)

$\theta$ is such that $\cos\theta = \frac{x}{\sqrt{x^2 + y^2}}$, $\sin\theta = \frac{y}{\sqrt{x^2 + y^2}}$, $\theta \in [0, 2\pi[$ (5.5)

$\varphi \in \left[-\frac{\pi}{2}, +\frac{\pi}{2}\right]$ such that $\varphi = \arctan\left(\frac{z}{\sqrt{x^2 + y^2}}\right)$ or $\varphi = \arcsin\left(\frac{z}{\sqrt{x^2 + y^2 + z^2}}\right)$

Same rules for $r'$, $\theta'$ and $\varphi'$ we just replace the unprimed symbols with the primed ones.

A *field* is a commutative ring (*F*, +, ·) such that 0 does not equal 1 and all elements of *F* except 0 have a multiplicative inverse.[10]

For ($\mathbb{H}$, +, ·) to be a field, the following field axioms must hold:

---

[9] Field (mathematics), Wikipedia, the free encyclopedia,
http://en.wikipedia.org/wiki/Field_(mathematics) Accessed 25 November 2003
[10] Ibid.

- ***Closure of $\mathbb{H}$ under + and ·***

For all $h$, $h'$ belonging to $\mathbb{H}$, both $h+h'$ and $h \cdot h'$ belong to $\mathbb{H}$ (or more formally, + and · are binary operations on $\mathbb{H}$);

$h = (x,y,z)$ and $h' = (x',y',z')$ \hfill (5.6)

**Operation +:** $h + h' = (x,y,z) + (x',y',z') = (x+x', y+y', z+z') = h''$ \hfill (5.7)

$h'' \in \mathbb{H}$

**Operation ·:** $h \cdot h' = (x,y,z) \cdot (x',y',z') = (r,\theta,\varphi) \cdot (r',\theta',\varphi')$

$$= re^{i\theta+j\varphi} \cdot r'e^{i\theta'+j\varphi'} = rr'e^{i(\theta+\theta')+j(\varphi+\varphi')}$$

$$= (rr', \theta+\theta', \varphi+\varphi') = h'' \quad (5.8)$$

If $\varphi + \varphi' > \dfrac{\pi}{2}$ or $\varphi + \varphi' < -\dfrac{\pi}{2}$ then $h'' = h \cdot h' = (rr', \theta+\theta'+\pi, \pi-\varphi-\varphi')$ \hfill (5.9)

$h'' \in \mathbb{H}$

Since $x+x', y+y', z+z'$, $rr'$, $\theta+\theta', \varphi+\varphi'$ are all real numbers, we see that $h+h'$ and $h \cdot h'$ are both members of $\mathbb{H}$. Thus $\mathbb{H}$ is closed under addition and multiplication.

- ***Both + and · are associative***

For all $h, h', h''$ in $\mathbb{H}$, $h + (h' + h'') = (h + h') + h''$ and $h \cdot (h' \cdot h'') = (h \cdot h') \cdot h''$ \hfill (5.10)

**Operation +:**

$h + (h' + h'') = (x,y,z) + ((x',y',z') + (x'',y'',z'')) = (x+(x'+x''), y+(y'+y''), z+(z'+z''))$

Since + is associative on $\mathbb{R}$:

$h + (h' + h'') = ((x+x')+x'', (y+y')+y'', (z+z')+z'')$ \hfill (5.11)

$h + (h' + h'') = (h + h') + h''$ \hfill (5.12)

+ is then associative on $\mathbb{H}$.

**Operation ·:**

$h \cdot (h' \cdot h'') = (x,y,z) \cdot [(x',y',z') \cdot (x'',y'',z'')] = (r,\theta,\varphi) \cdot [(r',\theta',\varphi') \cdot (r'',\theta'',\varphi'')]$ \hfill (5.13)

$h \cdot (h' \cdot h'') = re^{i\theta+j\varphi} \cdot (r'e^{i\theta'+j\varphi'} \cdot r''e^{i\theta''+j\varphi''}) = r(r'r'')e^{i[\theta+(\theta'+\theta'')]+j[\varphi+(\varphi'+\varphi'')]}$ \hfill (5.14)

Since + and × are associative on $\mathbb{R}$:

$h \cdot (h' \cdot h'') = (rr')r''e^{i[(\theta+\theta')+\theta'']+j[(\varphi+\varphi')+\varphi'']}$ \hfill (5.15)

$h \cdot (h' \cdot h'') = [(r,\theta,\varphi) \cdot (r',\theta',\varphi')] \cdot (r'',\theta'',\varphi'')$ \hfill (5.16)

$h \cdot (h' \cdot h'') = [(x,y,z) \cdot (x',y',z')] \cdot (x'',y'',z'') = (h \cdot h') \cdot h''$ \hfill (5.17)

· is then associative on $\mathbb{H}$.

- ***Both + and · are commutative***

For all $h, h'$ belonging to $\mathbb{H}$, $h + h' = h' + h$ and $h \cdot h' = h' \cdot h$ \hfill (5.18)

**Operation +:**

$h + h' = (x,y,z) + (x',y',z') = (x+x', y+y', z+z')$ \hfill (5.19)

Since + is commutative on $\mathbb{R}$

$h + h' = (x'+x, y'+y, z'+z) = h' + h$ \hfill (5.20)

+ is then commutative on $\mathbb{H}$.

**Operation ·:**

$h \cdot h' = (x,y,z) \cdot (x',y',z') = (r,\theta,\varphi) \cdot (r',\theta',\varphi')$ \hfill (5.21)

$h \cdot h' = re^{i\theta+j\varphi} \cdot r'e^{i\theta'+j\varphi'} = rr'e^{i(\theta+\theta')+j(\varphi+\varphi')}$ \hfill (5.22)

Since + and × are commutative on $\mathbb{R}$.

$h \cdot h' = r're^{i(\theta'+\theta)+j(\varphi'+\varphi)} = (r',\theta',\varphi') \cdot (r,\theta,\varphi) = h' \cdot h$ \hfill (5.23)

Operation · is then commutative.

- ### *Existence of an additive identity*

There exists an element 0 in $\mathbb{H}$, such that for all $h$ belonging to $\mathbb{H}$, $h + 0 = 0+h = h$.
The hypercomplex number $0 = (0,0,0)$ serves as the additive identity element for $\mathbb{H}$ because if $h = (x, y, z)$ is any hypercomplex number then

$$h + 0 = (x, y, z) + (0,0,0) = (x+0, y+0, z+0) = (x, y, z) = h \qquad (5.24)$$

The additive identity is the element (0,0,0) in Cartesian coordinates; or (0, $\theta$, $\varphi$) in spherical coordinates such that $\theta$ and $\varphi$ are any reals (modulus $r = 0$ in spherical coordinates).

- ### *Existence of a multiplicative identity*

There exists an element 1 in $\mathbb{H}$ different from 0, such that for all $h$ belonging to $\mathbb{H}$,
$$h \cdot 1 = 1 \cdot h = h \qquad (5.25)$$

The hypercomplex number $1 = (1,0,0)$ in spherical coordinates serves as the multiplicative identity element for $\mathbb{H}$ because if $h = (r, \theta, \varphi)$ is any hypercomplex number then

$$h \cdot 1 = (r, \theta, \varphi) \cdot (1,0,0) = (r \times 1, \theta + 0, \varphi + 0) = (r, \theta, \varphi) = h \qquad (5.26)$$

The number $1 = (1,0,0)$ in Cartesian coordinates.

We see from the above that $0 \neq 1$. The requirement $0 \neq 1$ ensures that the set which only contains a single zero is not a field.

- ### *Existence of additive inverses*

For every $h$ belonging to $\mathbb{H}$, there exists an element $-h$ in $\mathbb{H}$, such that:
$$h + (-h) = (-h) + h = 0 \qquad (5.27)$$

$$h = (x, y, z), -h = (a, b, c) \qquad (5.28)$$
$$h + (-h) = (x+a, y+b, z+c) = (0,0,0) \qquad (5.29)$$

this means $x+a=0$ and $y+b=0$ and $z+c=0$
$$\Rightarrow a = -x, b = -y, c = -z \qquad (5.30)$$

Therefore, the additive inverse of $h=(x,y,z)$ is $-h=(-x, -y, -z)$ in Cartesian coordinates.
In spherical coordinates the additive inverse of $h = (r, \theta, \varphi)$ is $-h = (r, \theta, \varphi + \pi)$

- ### *Existence of multiplicative inverses*

For every $h \neq 0$ belonging to $\mathbb{H}$, there exists an element $h^{-1}$ in $\mathbb{H}$, such that
$$h \cdot h^{-1} = h^{-1} \cdot h = 1. \qquad (5.31)$$

$$h = (r, \theta, \varphi), \ h^{-1} = (r', \theta', \varphi') \qquad (5.32)$$

$$h \cdot h^{-1} = re^{i\theta + j\varphi} \cdot r'e^{i\theta' + j\varphi'} = rr'e^{i(\theta+\theta') + j(\varphi+\varphi')} = 1 \cdot e^{i.0 + j.0} = 1 \qquad (5.33)$$

This means: $rr' = 1, \theta + \theta' = 0, \varphi + \varphi' = 0 \Leftrightarrow r' = \dfrac{1}{r}, \theta' = -\theta, \varphi' = -\varphi \qquad (5.34)$

The multiplicative inverse of $(r, \theta, \varphi)$ is then $\left(\dfrac{1}{r}, -\theta, -\varphi\right)$.

In Cartesian coordinates:

$$h = (x, y, z) \; ; \; h = re^{i\theta + j\varphi} \; ; \; r = \sqrt{x^2 + y^2 + z^2} \tag{5.35}$$

$$h^{-1} = \frac{1}{h} = \frac{1}{r} e^{-i\theta - j\varphi} \tag{5.36}$$

$$x' = \frac{1}{r}\cos\theta\cos\varphi = \frac{r}{r^2}\cos\theta\cos\varphi = \frac{x}{r^2}$$

$$y' = -\frac{1}{r}\sin\theta\cos\varphi = -\frac{r}{r^2}\sin\theta\cos\varphi = \frac{-y}{r^2} \tag{5.37}$$

$$z' = -\frac{1}{r}\sin\varphi = -\frac{r}{r^2}\sin\varphi = \frac{-z}{r^2}$$

$$h^{-1} = (\frac{x}{x^2 + y^2 + z^2}, \frac{-y}{x^2 + y^2 + z^2}, \frac{-z}{x^2 + y^2 + z^2}) \tag{5.38}$$

$$h^{-1} = \text{Inv}(h) \tag{5.39}$$

Directly from the axioms, one may show that $(\mathbb{H}, +)$ and $(\mathbb{H} - \{0\}, \cdot)$ are commutative groups and that therefore the additive inverse $-h$ and the multiplicative inverse $h^{-1}$ are uniquely determined by $h$. Furthermore, the multiplicative inverse of a product is equal to the product of the inverses:

$$(a \cdot b)^{-1} = a^{-1} \cdot b^{-1} \tag{5.40}$$

provided both $a$ and $b$ are non-zero. Other useful rules include

$$-a = (-1) \cdot a \tag{5.41}$$

and more generally

$$-(a \cdot b) = (-a) \cdot b = a \cdot (-b) \tag{5.42}$$

as well as

$$a \cdot 0 = 0, \tag{5.43}$$

all rules familiar from elementary arithmetic.[11]

All field axioms are satisfied except for distributivity of · over + as shown next.
- ***Is the operation · distributive over the operation + ?***

For all $h, h', h''$ belonging to $\mathbb{H}$, $h \cdot (h' + h'') = (h \cdot h') + (h \cdot h'')$ ?

$(x_1, y_1) \neq (0, 0)$, $(x_2 + x_3, y_2 + y_3) \neq (0, 0)$

$$p = h \cdot (h' + h'') = (x, y, z) \cdot [(x', y', z') + (x'', y'', z'')] = (x, y, z) \cdot (x' + x'', y' + y'', z' + z'') \tag{5.44}$$

$$p = (x_p, y_p, z_p) \tag{5.45}$$

$$x_p = (x(x' + x'') - y(y' + y''))(1 - \frac{z(z' + z'')}{\sqrt{x^2 + y^2}\sqrt{(x' + x'')^2 + (y' + y'')^2}})$$

$$y_p = (x(y' + y'') - y(x' + x''))(1 - \frac{z(z' + z'')}{\sqrt{x^2 + y^2}\sqrt{(x' + x'')^2 + (y' + y'')^2}}) \tag{5.46}$$

$$z_p = \sqrt{x^2 + y^2}(z' + z'') + z\sqrt{(x' + x'')^2 + (y' + y'')^2}$$

$$p' = (h_1 \cdot h_2) + (h_1 \cdot h_3) = (x_1, y_1, z_1) \cdot (x_2, y_2, z_2) + (x_1, y_1, z_1) \cdot (x_3, y_3, z_3) \tag{5.47}$$

$$p' = (x_{p'}, y_{p'}, z_{p'}) \tag{5.48}$$

- Let's study the $z$-component of the products $p$ and $p'$:

---

[11] Ibid.

$$z_p = \sqrt{x^2 + y^2}(z' + z'') + z\sqrt{(x' + x'')^2 + (y' + y'')^2} \tag{5.49}$$

$$z_p = \rho(z' + z'') + z\sqrt{\rho'^2 + \rho''^2 + 2(x'x'' + y'y'')} \tag{5.50}$$

From $p' = h \cdot h' + h \cdot h''$ we get:

$$z_{p'} = \rho z' + \rho' z + \rho z'' + \rho'' z \tag{5.51}$$

$$\boxed{\begin{aligned} z_{p'} &= \rho(z' + z'') + (\rho' + \rho'')z \\ z_p &= \rho(z' + z'') + z\sqrt{\rho'^2 + \rho''^2 + 2(x'x'' + y'y'')} \end{aligned}} \tag{5.52}$$
$$\tag{5.53}$$

$z_p = z_{p'}$ if and only if :

$$2\rho'\rho'' = 2(x'x'' + y'y'') = 2\vec{\rho}'.\vec{\rho}'' = 2\rho'\rho''.\cos\Delta\theta \text{ with } \Delta\theta = \theta'' - \theta' \tag{5.54}$$

$$\Rightarrow \cos\Delta\theta = 1 \Rightarrow \Delta\theta = 0 + 2k\pi, \ k \in \mathbb{Z} \tag{5.55}$$

$$\Rightarrow \theta'' = \theta' + 2k\pi, \ k \in \mathbb{Z} \tag{5.56}$$

$h'$ and $h''$ must then have the same argument $\theta$.

Therefore multiplication $\cdot$ is generally not distributive over addition +.

- Let's study the $x$-component.

$$x_{p'} = (xx' - yy')(1 - \frac{zz'}{\sqrt{x^2 + y^2}\sqrt{x'^2 + y'^2}}) + (xx'' - yy'')(1 - \frac{zz''}{\sqrt{x^2 + y^2}\sqrt{x''^2 + y''^2}}) \tag{5.57}$$

$$x_{p'} = (xx' - yy')(1 - \frac{zz'}{\rho\rho'}) + (xx'' - yy'')(1 - \frac{zz''}{\rho\rho''}) \tag{5.58}$$

But from the first formula (5.46) we have:

$$x_p = (x(x' + x'') - y(y' + y''))\left(1 - \frac{z(z' + z'')}{\rho\sqrt{\rho'^2 + \rho''^2 + 2(x'x'' + y'y'')}}\right) \tag{5.59}$$

$$x_p = ((xx' - yy') + (xx'' - yy''))\left(1 - \frac{z(z' + z'')}{\rho\sqrt{\rho'^2 + \rho''^2 + 2(x'x'' + y'y'')}}\right) \tag{5.60}$$

Let's express $x_p - x_{p'}$:

$$\begin{aligned} x_p - x_{p'} = (xx' - yy')\left(\frac{zz'}{\rho\rho'} - \frac{z(z' + z'')}{\rho\sqrt{\rho'^2 + \rho''^2 + 2(x'x'' + y'y'')}}\right) + \\ (xx'' - yy'')\left(\frac{zz''}{\rho\rho''} - \frac{z(z' + z'')}{\rho\sqrt{\rho'^2 + \rho''^2 + 2(x'x'' + y'y'')}}\right) = 0 \end{aligned} \tag{5.61}$$

$$\begin{aligned} \Rightarrow (xx' - yy')\frac{z}{\rho}\left(\frac{z'}{\rho'} - \frac{(z' + z'')}{\sqrt{\rho'^2 + \rho''^2 + 2(x'x'' + y'y'')}}\right) + \\ (xx'' - yy'')\frac{z}{\rho}\left(\frac{z''}{\rho''} - \frac{(z' + z'')}{\sqrt{\rho'^2 + \rho''^2 + 2(x'x'' + y'y'')}}\right) = 0 \end{aligned} \tag{5.62}$$

If $(z, \rho)$ is different from $(0,0)$

$$\Rightarrow (xx' - yy')\left(\frac{z'}{\rho'} - \frac{(z'+z'')}{\sqrt{\rho'^2 + \rho''^2 + 2(x'x'' + y'y'')}}\right) +$$

$$(xx'' - yy''))\left(\frac{z''}{\rho''} - \frac{(z'+z'')}{\sqrt{\rho'^2 + \rho''^2 + 2(x'x'' + y'y'')}}\right) = 0 \quad (5.63)$$

$$\Rightarrow (xx' - yy')\frac{z'}{\rho'} + (xx'' - yy'')\frac{z''}{\rho''} = ((xx' - yy') + (xx'' - yy''))\left(\frac{(z'+z'')}{\sqrt{\rho'^2 + \rho''^2 + 2(x'x'' + y'y'')}}\right)$$

Whatever $(xx' - yy')$ and $(xx'' - yy'')$ be, the equation is verified when

$$\frac{z'}{\rho'} = \frac{z''}{\rho''} = \frac{z' + z''}{\sqrt{\rho'^2 + \rho''^2 + 2(x'x'' + y'y'')}} \quad (5.64)$$

$\tan(\varphi') = \tan(\varphi'')$

means $\varphi' = \varphi'' + k\pi, \ k \in \mathbb{Z}$ (5.65)

but from z-component's study we got $(x'x'' + y'y'') = \rho'\rho''$, this means

$$\frac{z'}{\rho'} = \frac{z''}{\rho''} = \frac{z' + z''}{\rho' + \rho''} \quad (5.66)$$

then $z'' = \dfrac{z'\rho''}{\rho'}$

$\dfrac{z'}{\rho'} = \dfrac{z' + z''}{\rho' + \rho''} \Rightarrow z'(1 + \dfrac{\rho''}{\rho'})\rho' = z'(\rho' + \rho'') \Rightarrow \rho' + \rho'' = \rho' + \rho''$ such that

$\rho' + \rho'' \neq 0, \rho' \neq 0, \rho'' \neq 0$

From (5.66): $z'(\rho' + \rho'') = \rho'(z' + z'') \Rightarrow z'\rho'' = \rho'z'' \Rightarrow \dfrac{z'}{\rho'} = \dfrac{z''}{\rho''}$

So for $x_p$ and $x_{p'}$ to be equal, the condition is that $h'$ and $h''$ must have the same $\theta$ and same $\varphi$ or $\varphi'' = \varphi' + k\pi, k \in \mathbb{Z}$ ($h'$ and $h''$ must be collinear with $O$, that's to say located on the same straight line passing by the origin $O$).

The same thing can be concluded from the study of y-components.

The 3D spherical hypercomplex multiplication is then generally not distributive over addition in $\mathbb{H}$.

$re^{i\theta + k\varphi}(r'e^{i\theta' + k\varphi'} + r''e^{i\theta'' + k\varphi''}) = re^{i\theta + k\varphi}r'e^{i\theta' + k\varphi'} + re^{i\theta + k\varphi}r''e^{i\theta'' + k\varphi''}$ is generally false

$(\mathbb{H}, +, \cdot)$ then is not a field.

However, from the aforementioned we see that $(\mathbb{H}, +)$ and $(\mathbb{H} - \{0\}, \cdot)$ are commutative (Abelian) groups.

$\mathbb{H}$ is an an Abelian group under addition, meaning:

1. $(a + b) + c = a + (b + c)$ (+is associative)
2. There is an element 0 in R such that $0 + a = a$ (0 is the **zero element**)
3. $a + b = b + a$ (+is commutative)

4. For each *a* in $\mathbb{H}$ there exists −*a* in $\mathbb{H}$ such that $a + (-a) = (-a) + a = 0$ (−*a* is the additive inverse of a)[12]

$\mathbb{H} - \{0\}$ is also an Abelian group under multiplication.

The same thing can be demonstrated with sets of N-dimensional hyperspherical hypercomplex numbers, $N > 3$.

## 6. 3D hypercomplex inverse calculation

$$h = (x,y,z) \; ; \; h = (r,\theta,\varphi) = re^{i\theta + j\varphi} \; ; \; r = \sqrt{x^2 + y^2 + z^2}$$

$$Inv(h) = \frac{1}{h} = \frac{1}{r}e^{-i\theta - j\varphi} = \left(\frac{1}{r}, -\theta, -\varphi\right) = (x', y', z') \tag{6.1}$$

$$x' = \frac{1}{r}\cos\theta\cos\varphi = \frac{x}{r^2}$$

$$y' = -\frac{1}{r}\sin\theta\cos\varphi = \frac{-y}{r^2} \tag{6.2}$$

$$z' = -\frac{1}{r}\sin\varphi = \frac{-z}{r^2}$$

$$Inv(h) = \left(\frac{x}{x^2 + y^2 + z^2}, \frac{-y}{x^2 + y^2 + z^2}, \frac{-z}{x^2 + y^2 + z^2}\right) \tag{6.3}$$

**Rule:** To invert a hypercomplex number, we invert the modulus and we take the opposites of all its arguments

## 7. 3D hypercomplex conjugates

We call the *3D hypercomplex conjugate* of the 3D hypercomplex number $h(r,\theta,\varphi)$ the number $\bar{h}(r',\theta',\varphi')$ such that:

$$p = h\bar{h} = \bar{h}h = r^2 \tag{7.1}$$

$$p = re^{i\theta + j\varphi} \cdot r'e^{i\theta' + j\varphi'} = rr'e^{i(\theta+\theta') + j(\varphi+\varphi')} = r^2 e^{i\cdot 0 + j\cdot 0} \tag{7.2}$$

This means $r' = r$, $\theta' = -\theta$, $\varphi' = -\varphi$. \hfill (7.3)

$$\bar{h} = x - iy - jz \tag{7.4}$$

$$h + \bar{h} = 2x \tag{7.5}$$

$$h - \bar{h} = 2(iy + jz) \tag{7.6}$$

Geometric interpretation: The hypercomplex conjugate is found by reflecting *h* across the real axis.

We call the *second hypercomplex conjugate* of h the number $\bar{h}_2$ such as:

$$p = h\bar{h}_2 = r^2 e^{j2\varphi} = \frac{h^2}{e^{i2\theta}} \tag{7.7}$$

---

[12] Abelian group, Wikipedia, the free encyclopedia, http://en.wikipedia.org/wiki/Abelian_group Accessed 25 Nov 2003

Hence $\bar{h}_2 = re^{-i\theta+j\varphi} = x - iy + jz$ (7.8)

$$\frac{h+\bar{h}_2}{2} = x + jz$$ (7.9)

$$\frac{h-\bar{h}_2}{2i} = y$$ (7.10)

Geometric interpretation: The second complex conjugate is found by reflecting $h$ across plane ($xOz$).

We call the *third hypercomplex conjugate* of h the number $\bar{h}_3$ such as:

$$p = h\bar{h}_3 = r^2 e^{i2\theta} = \frac{r^2 \cos^2\varphi \cdot e^{i2\theta}}{\cos^2\varphi} = \frac{(x+jy)^2}{\cos^2\varphi}$$ (7.11)

Hence $\bar{h}_3 = re^{i\theta-j\varphi} = x + iy - jz$ (7.12)

$$\frac{h+\bar{h}_3}{2} = x + iy = r\cos\varphi \cdot e^{i\theta} \text{ is the complex part of } h.$$

$$\frac{h-\bar{h}_3}{2j} = z = r\sin\varphi \text{ is the second imaginary part}$$

Geometric interpretation: The third complex conjugate is found by reflecting $h$ across the complex plane ($xOy$).

## 8. 3D spherical hypercomplex number roots

We use the geometric form for the 3D spherical hypercomplex number root calculation. The $n^{th}$ hypercomplex root $h'(r',\theta',\varphi')$ of the hypercomplex number $h(r,\theta,\varphi)$ is defined by:

$h'^n = h$
(8.1)

$$\left(r'e^{i\theta'+j\varphi'}\right)^n = r'^n e^{in\theta'+jn\varphi'} = re^{i\theta+j\varphi}$$ (8.2)

$$\Rightarrow r'^n = r, \; n\theta' = \theta + 2k\pi, \; n\varphi' = \varphi + 2k\pi$$ (8.3)

$$\Rightarrow r' = \sqrt[n]{r}, \theta' = \frac{\theta}{n} + \frac{2k}{n}\pi, \varphi' = \frac{\varphi}{n} + \frac{2k}{n}\pi, k = 0,...,n-1$$ (8.4)

Hence we have $n$ different values of $\theta'$ and also $n$ different values of $\varphi'$. We have then $n \times n$ possible combinations of the two arguments.
We must not forget that $h(r,\theta,\varphi)$ has a replicate $(r, \theta+\pi, \pi-\varphi)$.
Therefore:

$$r' = \sqrt[n]{r}, \theta' = \frac{\theta+\pi}{n} + \frac{2k}{n}\pi, \varphi' = \frac{\pi-\varphi}{n} + \frac{2k}{n}\pi$$ (8.5)

$$\theta' = \frac{\theta}{n} + \frac{(2k+1)\pi}{n}, \varphi' = -\frac{\varphi}{n} + \frac{(2k+1)\pi}{n}$$ (8.6)

Hence we get *n* other different values of $\theta'$ and also *n* other different values of $\varphi'$. We have then $n \times n$ other possible combinations of the two arguments.

The roots resulting from the ordinary and replicate form may not be distinct. They can be double.

As a conclusion any 3D spherical hypercomplex number has at least *n* distinct 3D hypercomplex roots and at most $2n^2$ 3D hypercomplex roots.

**Example :** Calculating the square 3D hypercomplex roots of a complex number $c(r, \theta, \varphi = 0)$.

$$n = 2 \Rightarrow r' = \sqrt{r}, \theta' = \frac{\theta}{2} + k\pi, \varphi' = k\pi, k = 0,1 \tag{8.7}$$

| $\varphi'$ | $\theta'$ | Root |
|---|---|---|
| 0 | $\frac{\theta}{2}$ | $(\sqrt{r}, \frac{\theta}{2})$ (complex root) |
| 0 | $\frac{\theta}{2} + \pi$ | $(\sqrt{r}, \frac{\theta}{2} + \pi)$ complex root |
| $\pi$ | $\frac{\theta}{2}$ | $(\sqrt{r}, \frac{\theta}{2}, \pi) \equiv (\sqrt{r}, \frac{\theta}{2} + \pi, \pi - \pi) = (\sqrt{r}, \frac{\theta}{2} + \pi)$ |
| $\pi$ | $\frac{\theta}{2} + \pi$ | $(\sqrt{r}, \frac{\theta}{2} + \pi, \pi) \equiv (\sqrt{r}, \frac{\theta}{2} + 2\pi, \pi - \pi) = (\sqrt{r}, \frac{\theta}{2})$ |

$\equiv$ sign means "equivalent to replicate".

For $\varphi = 0$ we get the usual 2 complex roots. For $\varphi = \pi$ we get also the usual 2 complex roots. The roots then are double.

Roots of replicate form:

$$\theta' = \frac{\theta}{2} + \frac{(2k+1)\pi}{2}, \varphi' = \frac{(2k+1)\pi}{2} \text{ with } k = 0,1 \tag{8.9}$$

| $\varphi'$ | $\theta'$ | Hypercomplex root |
|---|---|---|
| $\frac{\pi}{2}$ | $\frac{\theta}{2} + \frac{\pi}{2}$ | $(\sqrt{r}, \frac{\theta}{2} + \frac{\pi}{2}, \frac{\pi}{2})$ |
| $\frac{\pi}{2}$ | $\frac{\theta}{2} + \frac{3\pi}{2}$ | $(\sqrt{r}, \frac{\theta}{2} + \frac{3\pi}{2}, \frac{\pi}{2})$ |
| $\frac{3\pi}{2} \equiv -\frac{\pi}{2}$ | $\frac{\theta}{2} + \frac{\pi}{2}$ | $(\sqrt{r}, \frac{\theta}{2} + \frac{\pi}{2}, -\frac{\pi}{2})$ |
| $\frac{3\pi}{2} \equiv -\frac{\pi}{2}$ | $\frac{\theta}{2} + \frac{3\pi}{2}$ | $(\sqrt{r}, \frac{\theta}{2} + \frac{3\pi}{2}, -\frac{\pi}{2})$ |

Roots 1 and 2 are replicates. Roots 3 and 4 are replicates too.

In Cartesian coordinates these 4 numbers are 2 opposite pure *z*-axis imaginary hypercomplex numbers.

## 9. Multivalued multiplication

Let's consider 3D spherical hypercomplex numbers. Because of spherical coordinates' duplicity (duplication/replication), a number *h(r,θ,φ)* can be similarly expressed: *h(r,θ,φ)* or *(r,θ+π,π–φ)*.
Let's consider two 3D hypercomplex numbers *h′* and *h″*.
We just forget all about the modulus, we consider both moduli for example to be equal to 1.
There are four cases for the product $h' \cdot h''$ depending on the replicates of *h′* and *h″* that we use. We just study $(\theta, \varphi)$ of the product $(\theta', \varphi') \cdot (\theta'', \varphi'')$.

|  | *h′* | *h″* | Product $h' \cdot h''$ |
|---|---|---|---|
| 1st case | $(\theta', \varphi')$ | $(\theta'', \varphi'')$ | $(\theta' + \theta'', \varphi' + \varphi'')$ |
| 2nd case | $(\theta', \varphi')$ | $(\theta'' + \pi, \pi - \varphi'')$ | $(\theta' + \theta'' + \pi, \pi - (\varphi'' - \varphi'))$ $\equiv (\theta' + \theta'', \varphi'' - \varphi')$ |
| 3rd case | $(\theta' + \pi, \pi - \varphi')$ | $(\theta'', \varphi'')$ | $(\theta' + \theta'' + \pi, \varphi'' - \varphi' + \pi)$ $\equiv (\theta' + \theta'', \varphi' - \varphi'')$ |
| 4th case | $(\theta' + \pi, \pi - \varphi')$ | $(\theta'' + \pi, \pi - \varphi'')$ | $(\theta' + \theta'', -\varphi' - \varphi'')$ |

We see then that the product is multivalued. It has the same longitude $\theta' + \theta''$, but the latitudes differ (4 latitudes / opposite 2 by 2) :
1st and 4th case : opposite latitudes
2nd and 3rd case : opposite latitudes

In the present theory we considered the first case only, that's to say $-\frac{\pi}{2} \leq \varphi \leq +\frac{\pi}{2}$.

## 10. Value of *j²*

To determine its value, we should determine which normal complex plane contains *j*, that's to say we should determine this plane's longitude θ. *j* is then function of θ and we note it *j(θ)*.

$$(j(\theta))^2 = j^2(\theta) = \left(je^{i\theta}\right)^2 = \left(e^{i\theta}e^{j\frac{\pi}{2}}\right)^2 = e^{i2\theta}e^{j\pi} = e^{i2\theta}(\cos\pi + j\sin\pi) = -e^{i2\theta} = e^{i(2\theta+\pi)}$$

Hence *j²* can be +1 or –1 or *i* or –*i* or whatever unit-modulus complex number depending on *j*'s longitude.

Also when two purely *j*-imaginary numbers are multiplied or a purely *j*-imaginary number is multiplied with a hypercomplex number, we must set the longitude for each of the *j*-imaginary numbers.

## 11. *N*-dimensional hyperspherical hypercomplex numbers

## 11.1. Hyperspherical Coordinates

Spherical coordinates are a generalization of polar coordinates, and can be further generalized to the *N*-sphere (or *N*-hypersphere) with (*N*–2) latitudes $\varphi_i$ and one longitude $\theta$ (or say *(N–1)* arguments of the hypercomplex number *h*). We call these coordinates in such a case *hyperspherical* coordinates.

Let's make $\theta = \theta_2, \theta_i = \varphi_i, i \geq 3, i$ natural

An *N*-dimensional hyperspherical hypercomplex number *h* is defined by the modulus r and *N*–1 arguments
$$h = (r, \theta_2, ..., \theta_N) \tag{11.1}$$

*h*'s Cartesian coordinates are defined by:
$$x_1 = r \cos\theta_2 \cos\theta_3 .... \cos\theta_N$$
$$x_2 = r \sin\theta_2 \cos\theta_3 ... \cos\theta_N$$
$$...$$
$$x_k = r \sin\theta_k \cos\theta_{k+1} .... \cos\theta_N \tag{11.2}$$
$$....$$
$$x_N = r \sin\theta_N$$

Where $r = r_N = \sqrt{\sum_{n=1}^{N} x_n^2}$ (11.3)

And the *n*-dimensional modulus $r_n = \sqrt{\sum_{k=1}^{n} x_k^2}$ (11.4)

We call $r_N = r$ the *N*-dimensional modulus or simply *h*'s modulus.
See Fig.11.1.

**<u>Important note</u>**
An *N*-dimensional hypercomplex number *h* has *(N–1)* arguments $\theta_2, ..., \theta_N$.
If we consider $\theta_2$ to be in the interval $I = \left[-\frac{\pi}{2}, +\frac{\pi}{2}\right]$, then we need to consider one extra argument $\theta_1$ which will have 2 values :
$+\frac{\pi}{2}$ for positive reals.
$-\frac{\pi}{2}$ for negative reals.
For zero it does not matter. (*x*-axis value is $x = r_1 \sin\theta_1$)
(in such a case $\sin\theta_2$ gives the 2D imaginary part; $\cos\theta_2$ gives the *x*-axis modulus, that is to say $r_1$, $\sin\theta_1$ gives the sign of the *x*-component)

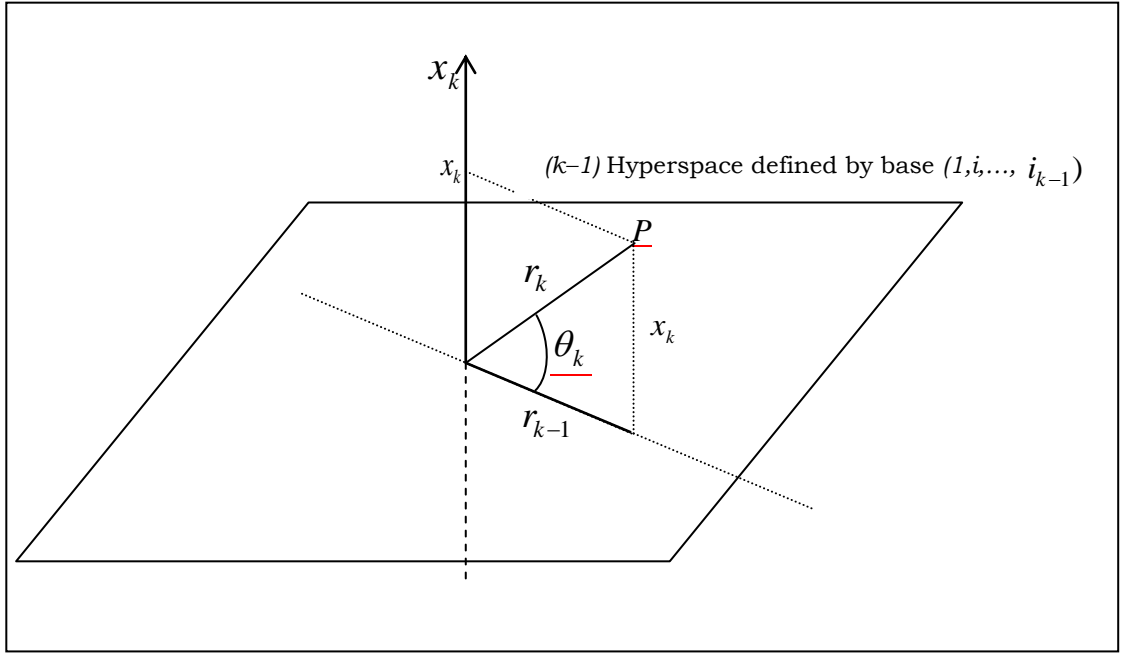

Figure 11.1. *Hyperspherical coordinates*

## 11.2. Four-Dimensional Hyperspherical Hypercomplex numbers

A 4D hypercomplex number $h$ is defined by its modulus $r$ and three arguments $\theta$, $\varphi$, $\psi$ (there will be an extra latitude $\psi$ to the 3D case).

$$h = re^{i\theta + j\varphi + k\psi} \tag{11.5}$$

$$h = re^{i\theta + j\varphi}(\cos\psi + k\sin\psi) \tag{11.6}$$

$$h = re^{i\theta + j\varphi}\cos\psi + kr\sin\psi \text{ because } e^{i\theta + j\varphi}k = k \tag{11.7}$$

$$h = r\cos\psi . e^{i\theta}(\cos\varphi + j\sin\varphi) + kr\sin\psi \tag{11.8}$$

$$h = r\cos\psi . e^{i\theta}\cos\varphi + jr\cos\psi\sin\varphi + kr\sin\psi \tag{11.9}$$

$$h = r\cos\psi\cos\varphi\cos\theta + ir\cos\psi\cos\varphi\sin\theta + jr\cos\psi\sin\varphi + kr\sin\psi \tag{11.10}$$

$h = (x, y, z, w) = (r, \theta, \varphi, \psi)$ and $h' = (x', y', z', w') = (r', \theta', \varphi', \psi')$ are two 4D hyperspherical hypercomplex numbers. $h'' = (x'', y'', z'', w'') = (r'', \theta'', \varphi'', \psi'')$ is their product.

$$h'' = h \cdot h' = (r'', \theta'', \varphi'', \psi'') = (rr', \theta + \theta', \varphi + \varphi', \psi + \psi') \tag{11.11}$$

As in constructing the 3D hypercomplex numbers, the 4D hypercomplex numbers are constructed through the orthogonal complex plane construction.

We take the 3D moduli $r_3 = \sqrt{x^2 + y^2 + z^2}$ and $r_3' = \sqrt{x'^2 + y'^2 + z'^2}$ of the 3D hypercomplex parts and we apply complex multiplication on them along with $w$ and $w'$.

We get $r_3 r_3' - ww'$ for 3 D space component. And $r_3 w' + r_3' w = w''$ for $w$-component.
We multiply then the 3D space component :
by $\cos(\varphi + \varphi')\cos(\theta + \theta')$ to obtain $x''$.
by $\cos(\varphi + \varphi')\sin(\theta + \theta')$ to obtain $y''$.
and by $\sin(\varphi + \varphi')$ to obtain $z''$.

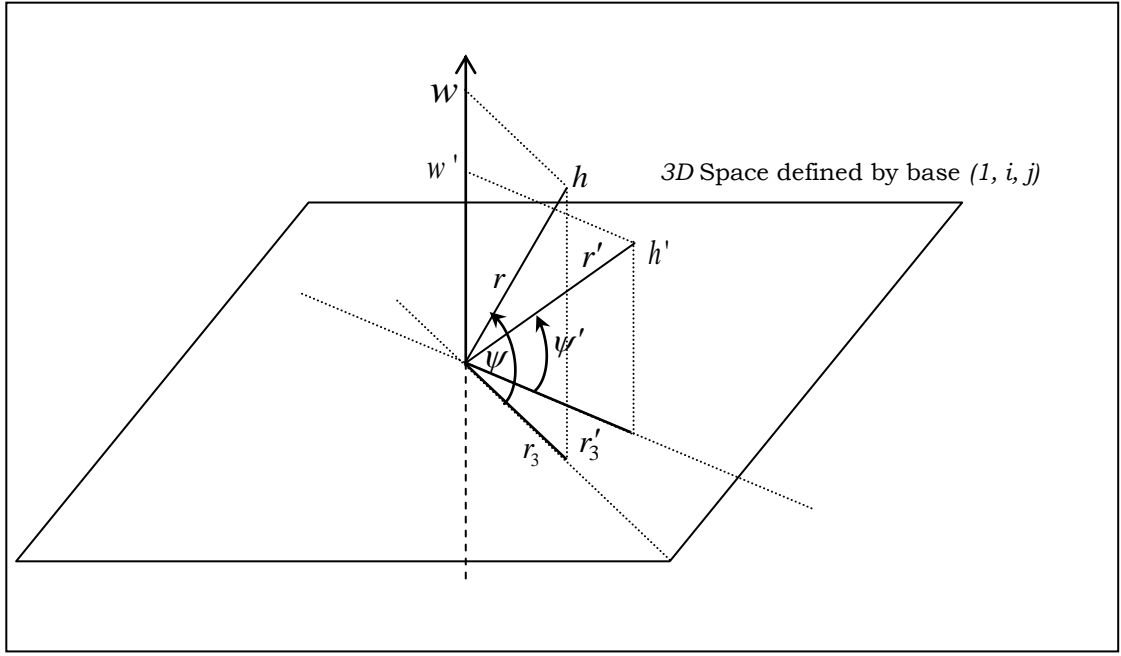

Fig.11.2. *Four dimensional hypercomplex multiplication*

## 11.3. *N*-dimensional generalized Euler's formula

We construct it using the complex-hypercomplex orthogonal construction. For any hypercomplex space *S(n)*, we can extend such a hyperspace using another plane $\pi(n+1)$ which contains axis $(Ox_{n+1})$ and the point *P*. Since axis $(Ox_{n+1})$ is normal to *n*-D hyperspace *S(n)*, the plane containing it and point *P* is thus normal to hyperspace *S(n)*.

$i_2,...,i_N$ are the unit *n*-imaginary numbers.
$i_2 = i$
A unit-modulus *N*-D hypercomplex number *h* can be expressed:

$$h = e^{\sum_{k=2}^{N} i_k \theta_k} = e^{i_2\theta_2 +...+ i_N\theta_N} = e^{i_2\theta_2}.e^{i_3\theta_3}....e^{i_N\theta_N} = \prod_{k=2}^{N} e^{i_k \theta_k} \quad (11.12)$$

$$e^{\sum_{k=2}^{N} i_k \theta_k} = \prod_{k=2}^{N} (\cos\theta_k + i_k \sin\theta_k), \ k \geq 2, k \in N \quad (11.13)$$

Equation (11.13) is the *N*-D generalized Euler's formula.
As in section 3.2.3.b., after developing equation (11.13) and comparing with equations (11.2) we can deduce the following theorem and fundamental property:

**General Theorem:** $i_n$ multiplied by an *(n–1)*-dimensional hypercomplex number of modulus 1 remains unchanged (or multiplied by a hypercomplex number located on the *(n–1)*-unit-hypersphere). $i_n$ is insensitive to rotation in the *(n–1)*-hypercomplex space. $1 \leq n \leq N$

**Fundamental property**

$$e^{\sum_{k=2}^{m-1} i_k \theta_k} . i_m = e^{\sum_{k=2}^{m-1} i_k \theta_k} . e^{i_m \frac{\pi}{2}} = i_m, \forall m \geq 3, m \in N \qquad (11.14)$$

$h$'s or Euler's formula Cartesian components (coordinates) are:

$$\boxed{\begin{aligned}
x_1 &= \prod_{k=2}^{N} \cos \theta_k \\
x_2 &= \sin \theta_2 \prod_{k=3}^{N} \cos \theta_k \\
&\ldots \\
x_m &= \sin \theta_m \prod_{k=m+1}^{N} \cos \theta_k \\
&\ldots \\
x_N &= \sin \theta_N
\end{aligned}} \qquad (11.15)$$

In fact, we sum up arguments multiplied by their respective unit imaginary numbers till infinity. But for any $k>N$, $\theta_k = 0$.

Hence we can rewrite:

$$e^{\sum_{k=2}^{N} i_k \theta_k} = e^{\sum_{k=2}^{\infty} i_k \theta_k} \qquad (11.16)$$

Hence

$$\prod_{k=N+1}^{\infty} \cos \theta_k = 1 \qquad (11.17)$$

This is the term which multiplies $i_N \sin \theta_N$ (which is 1).

**Important note**

If we consider $\theta_1$ having the values $\frac{\pi}{2}$ for strictly positive numbers and $-\frac{\pi}{2}$ for strictly negative numbers, for zero it doesn't matter. $x$-axis values will be generated using its modulus (absolute value) and the argument $\theta_1$ (this axis will represent imaginary numbers with respect to another 0-dimension axis).

**In this case** $\theta_2$ must be comprised between $-\frac{\pi}{2}$ and $\frac{\pi}{2}$.

## 11.4. Multiplication: *N*-dimensional hypercomplex general case

An *N*-D hypercomplex number is defined by its modulus $r$ and its *(N–1)* arguments $\theta_2, \ldots, \theta_N$.

The unit imaginary numbers $i_k$, $k=2,\ldots,N$ have the following interesting property:
Every unit *m*-imaginary number multiplied by a unit-modulus (*m*–1)-number remains invariant.

$$e^{\sum_{n=2}^{m-1} i_n \theta_n} i_m = i_m e^{\sum_{n=2}^{m-1} i_n \theta_n} = i_m \qquad (11.18)$$

Note: Case of imaginary complex $i$: $1.i = i.1 = i$

$$h = r \prod_{n=2}^{N} e^{i_n \theta_n}, \quad h' = r' \prod_{n=2}^{N} e^{i_n \theta'_n}, \quad (11.19)$$

$$\text{Product} = h'' = hh' = rr' \prod_{n=2}^{N} e^{i_n (\theta_n + \theta'_n)} \quad (11.20)$$

Hence multiplying $N$-D hyperspherical hypercomplex numbers results in:
- multiplying moduli
- $i^{\text{th}}$-order argument = sum of $i^{\text{th}}$-order arguments.

This means there occurs a rotation of the first hypercomplex $h$ (vector $\vec{v}$) in the $i^{\text{th}}$-order plane by an angle equal to the $i^{\text{th}}$-order number argument.

Hence hypercomplex multiplication is very useful to express rotations about a point in any $N$-D hyperspace.

$h$'s Cartesian coordinates in function of hyperspherical coordinates $(r, \theta_2, ..., \theta_N)$:

$$x = x_1 = r \prod_{n=2}^{N} \cos \theta_n$$

$$y = x_2 = r \sin \theta_2 \prod_{n=3}^{N} \cos \theta_n$$

$$z = x_3 = r \sin \theta_3 \prod_{n=4}^{N} \cos \theta_n \quad (11.21)$$

....

$$x_{N-1} = r \sin \theta_{N-1} \cos \theta_N$$

$$x_N = r \sin \theta_N$$

Where $r = r_N = \sqrt{\sum_{n=1}^{N} x_n^2} \quad (11.22)$

And the $n$-dimensional modulus $r_n = \sqrt{\sum_{k=1}^{n} x_k^2} \quad (11.23)$

Cartesian coordinates (components) of the product are:
$h''(x''_1, ..., x''_N)$

$$x''_1 = rr' \prod_{n=2}^{N} \cos(\theta_n + \theta'_n)$$

$$x''_2 = rr' \sin(\theta_2 + \theta'_2) \prod_{n=3}^{N} \cos(\theta_n + \theta'_n) \quad (11.24)$$

....

$$x''_{N-1} = rr' \sin(\theta_{N-1} + \theta'_{N-1}) \cos(\theta_N + \theta'_N)$$

$$x''_N = rr' \sin(\theta_N + \theta'_N)$$

## 11.5. Expression of the product in function of the Cartesian coordinates

$h(x_1,...,x_N) = h(r, \theta_2,...,\theta_N)$, $h'(x'_1,...,x'_N) = h'(r', \theta'_2,...,\theta'_N)$ are two $N$-dimensional hypercomplex numbers.

*__Note: Moduli's fundamental property__*

$r = r_N$

$r_{N-1} = r_N \cos\theta_N = r\cos\theta_N$

$r_{N-2} = r_{N-1} \cos\theta_{N-1} = r\cos\theta_N \cos\theta_{N-1}$

... (11.25)

$r_k = r\prod_{n=1}^{k} \cos\theta_{N-k} = r\prod_{n=k+1}^{N} \cos\theta_n$

...

$r_1 = x = r \cdot$ product of cosines of all arguments

Let's rewrite the terms $\cos(\theta_n + \theta'_n)$ and $\sin(\theta_n + \theta'_n)$ in (11.24) in function of trigonometric terms $\cos\theta_n$, $\cos\theta'_n$, $\sin\theta_n$ and $\sin\theta'_n$

$$\cos(\theta_n + \theta'_n) = \cos\theta_n \cos\theta'_n - \sin\theta_n \sin\theta'_n = \cos\theta_n \cos\theta'_n (1 - \frac{\sin\theta_n \sin\theta'_n}{\cos\theta_n \cos\theta'_n})$$

$$= \cos\theta_n \cos\theta'_n (1 - \tan\theta_n \tan\theta'_n) \quad (11.26)$$

$$\sin(\theta_n + \theta'_n) = \sin\theta_n \cos\theta'_n + \cos\theta_n \sin\theta'_n, \quad \forall n \geq 2 \quad (11.27)$$

Where $\tan\theta_n = \dfrac{x_n}{r_{n-1}}$, and $\tan\theta'_n = \dfrac{x'_n}{r'_{n-1}}$, $n \geq 2$ (see figure 11.1.)

For the special case n=2 (conventional complex plane), we consider $r_1 = x_1 = x$ hence $\tan\theta_2 = \dfrac{y}{x}$ and $\tan\theta'_2 = \dfrac{y'}{x'}$

We can rewrite equations (11.24) as follows:

$$x''_1 = rr'\prod_{n=2}^{N}[\cos\theta_n \cos\theta'_n (1 - \tan\theta_n \tan\theta'_n)]$$

$$x''_2 = rr'(\sin\theta_2 \cos\theta'_2 + \cos\theta_2 \sin\theta'_2)\prod_{n=3}^{N}[\cos\theta_n \cos\theta'_n (1 - \tan\theta_n \tan\theta'_n)]$$

.....

$$x''_k = rr'(\sin\theta_k \cos\theta'_k + \cos\theta_k \sin\theta'_k)\prod_{n=k+1}^{N}[\cos\theta_n \cos\theta'_n (1 - \tan\theta_n \tan\theta'_n)] \quad (11.28)$$

.....

$x''_{N-1} = rr'(\sin\theta_{N-1} \cos\theta'_{N-1} + \cos\theta_{N-1} \sin\theta'_{N-1})\cos\theta_N \cos\theta'_N (1 - \tan\theta_N \tan\theta'_N)$

$x''_N = rr'(\sin\theta_N \cos\theta'_N + \cos\theta'_N \sin\theta'_N)$

Knowing that:

$h$'s coordinates are:

$$x = x_1 = r\prod_{n=2}^{N} \cos\theta_n$$

$$y = x_2 = r \sin \theta_2 \prod_{n=3}^{N} \cos \theta_n$$

$$z = x_3 = r \sin \theta_3 \prod_{n=4}^{N} \cos \theta_n$$

.... (11.29)

$$x_n = r \sin \theta_n \prod_{k=n}^{N} \cos \theta_n$$

....

$$x_{N-1} = r \sin \theta_{N-1} \cos \theta_N$$

$$x_N = r \sin \theta_N$$

The same thing for $h'$ adding only the "prime" symbol.

We can rewrite the Cartesian coordinates of the product:

$$x_1'' = x'' = r \prod_{n=2}^{N} \cos \theta_n \cdot r' \prod_{n=2}^{N} \cos \theta_n' \cdot \prod_{n=2}^{N} (1 - \tan \theta_n \tan \theta_n')$$

$$x_1'' = x_1 x_1' \prod_{n=2}^{N} (1 - \tan \theta_n \tan \theta_n') = xx' \prod_{n=2}^{N} (1 - \frac{x_n x_n'}{r_{n-1} r_{n-1}'}) = (xx' - yy') \prod_{n=3}^{N} (1 - \frac{x_n x_n'}{r_{n-1} r_{n-1}'})$$

$$x_2'' = (xy' + yx') \prod_{n=3}^{N} (1 - \tan \theta_n \tan \theta_n') = (x_1 x_2' + x_2 x_1') \prod_{n=3}^{N} (1 - \frac{x_n x_n'}{r_{n-1} r_{n-1}'})$$

.....

$$x_k'' = (x_k r_{k-1}' + x_k' r_{k-1}) \prod_{n=k+1}^{N} (1 - \tan \theta_n \tan \theta_n') = (x_k r_{k-1}' + x_k' r_{k-1}) \prod_{n=k+1}^{N} (1 - \frac{x_n x_n'}{r_{n-1} r_{n-1}'})$$

.....

$$x_{N-1}'' = rr'(\sin \theta_{N-1} \cos \theta_{N-1} + \cos \theta_{N-1}' \sin \theta_{N-1}') \cos \theta_N \cos \theta_N' (1 - \tan \theta_N \tan \theta_N')$$

$$x_{N-1}'' = (x_{N-1} r_{N-2}' + x_{N-1}' r_{N-2}) (1 - \frac{x_N x_N'}{r_{N-1} r_{N-1}'})$$ (11.30)

$$x_N'' = rr'(\sin \theta_N \cos \theta_N + \cos \theta_N' \sin \theta_N')$$

$$x_N'' = x_N r_{N-1}' + x_N' r_{N-1}$$

Hence:

$$\boxed{\begin{aligned}
x_1'' &= x'' = (x_1 x_1' - x_2 x_2') \prod_{n=3}^{N} (1 - \frac{x_n x_n'}{r_{n-1} r_{n-1}'}) = (xx' - yy') \prod_{n=3}^{N} (1 - \frac{x_n x_n'}{r_{n-1} r_{n-1}'}) = x_1 x_1' \prod_{n=2}^{N} (1 - \frac{x_n x_n'}{r_{n-1} r_{n-1}'}) \\
x_2'' &= y'' = (x_1 x_2' + x_2 x_1') \prod_{n=3}^{N} (1 - \frac{x_n x_n'}{r_{n-1} r_{n-1}'}) = (xy' + yx') \prod_{n=3}^{N} (1 - \frac{x_n x_n'}{r_{n-1} r_{n-1}'}) \\
&\ldots \\
x_k'' &= (x_k r_{k-1}' + x_k' r_{k-1}) \prod_{n=k+1}^{N} (1 - \frac{x_n x_n'}{r_{n-1} r_{n-1}'}) = (x_k r_{k-1}' + x_k' r_{k-1}) \prod_{n=k}^{N-1} (1 - \frac{x_{n+1} x_{n+1}'}{r_n r_n'}) \\
&\ldots \\
x_N'' &= x_N r_{N-1}' + x_N' r_{N-1}
\end{aligned}}$$ (11.31)

Above is the expression of the Cartesian coordinates of the product $h''$ of $h$ and $h'$ in function of the Cartesian coordinates of $h$ and $h'$.

### 11.5.1. Division by zero: Indetermination of Cartesian components

The $r_n = 0$ problem happens only if the $m$ first consecutive coordinates are equal to 0. That's to say if the $(m+1)$-hypercomplex part is purely $i_{m+1}$ imaginary.

**Proof**

$$r_n = \sqrt{\sum_{k=1}^{n} x_k^2}, \; n \le m \qquad (11.32)$$

$r_n = 0$ if and only if $x_k = 0, \forall k \le n$. End of proof.

In that case the $m$ first $r_n$'s are equal to zero. Starting from $m+1$ ($x_{n+1} \ne 0$), $r_k$ will never get equal to zero.

Hence there will be no indetermination problem for $n \ge m+1$ components.

Therefore, the indetermination problem will be posed for the first $m$ coordinates. The $m-1$ first arguments $(\theta_2, ..., \theta_n)$ must therefore be provided.

$\theta_{m+1}$ will necessarily be equal to $\pm \dfrac{\pi}{2}$. $h_{m+1}$ hypercomplex number must be a purely $(m+1)$-imaginary number.

## 11.6. Is ($\mathbb{H}_N$,+,·) a field?

($\mathbb{H}_N$,+,·) is the set $\mathbb{H}_N$ of *N*-dimensional hypercomplex numbers equipped with the above defined laws of addition + and multiplication ·. Is that set a field??
No, It isn't. Similarly to the 3D case, it can be demonstrated that multiplication generally does not distribute over addition. It distributes only for vectors which are collinear or parallel (hypercomplex numbers located on the same straight line passing by the origin).

## 11.7. Inverse general case: *N*-D hypercomplex numbers

$$h = (r, \theta_2, ..., \theta_N) = (x_1, x_2, ..., x_N) \qquad (11.33)$$

**Hyperspherical coordinates**

$$\text{Inv}(h) = \frac{1}{h} = h^{-1} = \frac{1}{re^{\sum_{k=2}^{N} i_k \theta_k}} = \frac{1}{r} e^{-\sum_{k=2}^{N} i_k \theta_k} = (\frac{1}{r}, -\theta_2, ..., -\theta_N) \qquad (11.34)$$

**Rule:** To invert an *N*-D hypercomplex number, we invert the modulus and we take the opposites of all its arguments

**Cartesian coordinates**

$$h = x_1 + \sum_{n=2}^{N} i_n x_n$$

$$\text{Inv}(h) = \frac{1}{h} = \frac{x_1}{r^2} + \sum_{n=2}^{N} i_n (\frac{-x_n}{r^2}) = \frac{x_1}{r^2} - \sum_{n=2}^{N} i_n \frac{x_n}{r^2} = (\frac{x_1}{r^2}, -\frac{x_2}{r^2}, ..., -\frac{x_N}{r^2}) \quad (11.35)$$

$$r = \sqrt{\sum_{n=1}^{N} x_n^2}$$

## 11.8. *N*-D Conjugate: same modulus, opposite arguments

We call the *N*-D hypercomplex conjugate of the *N*-D hyperspherical hypercomplex number $h(r, \theta_2, ..., \theta_N)$ the number $\bar{h}(r', \theta'_2, ..., \theta'_N)$ such that:

$$p = h\bar{h} = \bar{h}h = r^2 \quad (11.36)$$

$$p = re^{\sum_{k=2}^{N} i_k \theta_k} \cdot r'e^{\sum_{k=2}^{N} i_k \theta'_k} = rr'e^{\sum_{k=2}^{N} i_k (\theta_k + \theta'_k)} = r^2 \quad (11.37)$$

This means $r' = r$, $\theta'_k = -\theta_k, k = 2, ..., N$. $\quad (11.38)$

$$\bar{h} = x_1 + \sum_{k=2}^{N} i_k (-x_k) = x_1 - \sum_{k=2}^{N} i_k x_k \quad (11.39)$$

$$h + \bar{h} = 2x_1 \quad (11.40)$$

$$h - \bar{h} = 2\sum_{k=2}^{N} i_k x_k \quad (11.41)$$

Geometric interpretation: The hypercomplex conjugate is found by reflecting *h* across the real axis.

Other conjugates can be obtained through the combination of reversed arguments.

## 11.9. *N*-dimensional duplicity (replication) of arguments

When argument $\theta_k$ is replaced with $\pi - \theta_k$, then argument $\theta_{k-1}$ must be replaced with $\theta_{k-1} + \pi$ to keep the lower index Cartesian coordinates unchanged. The other lower index arguments don't need to be changed.

So when we calculate roots, we need to consider replicates defined by $\pi - \theta_k$, $k \geq 3$, and $\theta_{k-1} + \pi$.

## 11.10. *N*-D hypercomplex number's m<sup>th</sup> hypercomplex roots

We use the geometric form for the *N*D hyperspherical hypercomplex number root calculation. The $m^{th}$-*N*D hypercomplex root $h'(r', \theta'_2, ..., \theta'_N)$ of the *N*D hypercomplex number $h(r, \theta_2, ..., \theta_N)$ is defined by:

$$h'^m = h \quad (11.42)$$

$$\left( r'e^{\sum_{k=2}^{N} i_k \theta'_k} \right)^m = r'^m e^{\sum_{k=2}^{N} i_k m \theta'_k} = r e^{\sum_{k=2}^{N} i_k \theta_k} \tag{11.43}$$

$$\Rightarrow r'^m = r,\ m\theta'_k = \theta_k + 2j\pi \tag{11.44}$$

$$\Rightarrow r' = \sqrt[m]{r},\ \theta'_k = \frac{\theta_k}{m} + \frac{2j}{m}\pi,\ j = 0,\ldots,m-1;\ k = 2,\ldots,N \tag{11.45}$$

There are $N-1$ arguments, and m roots for every argument. This makes $m^{N-1}$ roots, in addition to the replicates' roots. Replicates are: $(r,\ldots,\theta_{k-1}+\pi,\pi-\theta_k,\ldots)$, $k=3,\ldots,N$.

Any $N$-dimensional hypercomplex number $h$ will have an infinity of $L$-D hypercomplex roots, $L \geq N$. We use the simple form and the replicate form.

Any complex number will have an infinity of $L$-D roots, $L \geq 2$. We use the simple form as well as the replicate one $(r,\ldots,\theta_{k-1}+\pi,\pi-\theta_k,\ldots)$, $k>2$, $\theta_k = 0$.

For the special case of square roots ($m=2$) of an $N$-D hypercomplex number, $\theta'_k = \frac{\pi}{2}, \frac{3\pi}{2}$, $k > N$ we will have 2 $i_k$-purely imaginary roots, for every $k > N$, one positive and another negative. For $N=2$ case of complex numbers: there are an infinity of positive and negative $i_k$-purely imaginary hypercomplex roots, $k > 2$.

## 11.11. Hypercomplex division in Cartesian coordinates

As in complex division, we multiply and divide by the denominator's conjugate. We get the modulus squared in the denominator, and we calculate the numerator using the product's Cartesian coordinates formula.

# 12. The new vector algebra: Merging number and vector theories in just one theory dealing with the same entity

In this section, we are going to deal with hypercomplex numbers as vectors.
The basis of an *N*-hypercomplex space is $(i_1 = 1, i_2 = i, ..., i_N)$.
In fact we can adopt vector notation and use the *N*-vector space orthonormed basis $(\vec{i}_1, \vec{i}_2, ..., \vec{i}_N)$

**3D case**

$$h = x + iy + jz \tag{12.1}$$

$$\vec{v} = x\vec{i} + y\vec{j} + z\vec{k} = re^{\vec{j}\theta + \vec{k}\varphi} \tag{12.2}$$

**N-D case**

$$h = \sum_{k=1}^{N} i_k x_k = re^{\sum_{k=2}^{N} i_k \theta_k} \tag{12.3}$$

$$\vec{v} = \sum_{k=1}^{N} \vec{i}_k x_k = re^{\sum_{k=2}^{N} \vec{i}_k \theta_k} \tag{12.4}$$

**Vector multiplication**
It is defined the same way as hypercomplex numbers.
For the first time in the history of mathematics, we can divide a vector by another.

The vector product thus defined is dependent on the referential unlike the usual dot and cross vector products which are independent.

**Scalar's Meaning**
In mathematics, the meaning of scalar depends on the context; it can refer to real numbers or complex numbers or rational numbers, or to members of some other specified field (mathematics). Generally, when a vector space over the field F is studied, then F is called the field of scalars.[13]

In physics a scalar is a quantity that can be described by a single number (either dimensionless, or in terms of some physical quantity). Scalar quantities have magnitude, but not a direction and should thus be distinguished from vectors. More formally, a scalar is a quantity that is invariant under coordinate rotations (or Lorentz transformations, for relativity). A scalar is formally a tensor of rank zero.[14]

Let's find the meaning of the usual scalars.
The celebrated formula for multiplying a vector $\vec{v}$ by a scalar *s* is:
$s \cdot \vec{v}(x, y, z, ..., x_N) = \vec{v}'(sx, sy, sz, ..., sx_N)$
Let's express this product in function of the new hypercomplex product.
We see that $\vec{v}$ does not change direction if *s* is positive. The resulting vector is opposite if *s* is negative.

---

[13] Scalar, Wikipedia, the free encyclopedia http://en.wikipedia.org/wiki/Scalar_(mathematics)
Accessed 09 May 2004
[14] Ibidem

Let's take the case of 1 and −1.
The scalar 1 will be represented by a vector whose modulus is 1 and all its $N-1$ arguments are 0.
$s = 1 = (1,0,0,...,\theta_N = 0)$ hyperspherical $= (1,0,0,...,0)$ Cartesian

$s=1 \equiv$ *unit vector of x-axis.*

When $s = -1$, $\vec{v}$ is inverted: $\vec{v}'(-x,-y,-z,...,-x_N)$.

## 1st Solution

$$x1 = \prod_{k=2}^{N} \cos\theta_k$$

$$x2 = \sin\theta_k \prod_{k=3}^{N} \cos\theta_k$$

...

$$x_n = \sin\theta_n \prod_{k=n+1}^{N} \cos\theta_k \qquad (12.5)$$

...

$$x_{N-1} = \sin\theta_{N-1} \cos\theta_N$$

$$x_N = \sin\theta_N$$

Changing $\theta_N$ with $\theta_N + \pi$ makes $\sin(\theta_N + \pi) = -\sin\theta_N$.
In addition: $\cos(\theta_N + \pi) = -\cos\theta_N$ which will be multiplied by all the lower index terms. Thus adding $\pi$ only to $\theta_N$ will change the whole $\vec{v}$ from $(x,y,z,...,x_N)$ to $(-x,-y,-z,...,-x_N)$.

This means that the scalar $s = -1$ is an N-D vector of modulus 1, of $\theta_N = \pi$ and the rest of its arguments are equal to zero.
Scalar $-1 = (r = 1,0,0,...,\theta_N = \pi)$ hyperspherical.

Let's calculate the value of this vector in Cartesian coordinates.
All $x_i, i \geq 2$ are null except for $x_1$ where the cosine product is equal to $-1$.

Hence the scalar $s = -1$ can be written in Cartesian coordinates $(-1,0,0,...,0)$.

The calculation of hyperspherical coordinates departing from Cartesian $(-1,0,0,...,0)$ gives:

$\theta_2 = \pi, \theta_k = 0, k = 3,...,N$ because $\cos\theta_k = \dfrac{r_k}{r_{k-1}} = \dfrac{1}{1} = 1 \Rightarrow \theta_k = 0$ or because

$\tan\theta_k = \dfrac{x_k}{r_{k-1}} = \dfrac{0}{1} = 0 \Rightarrow \theta_k = 0$.

Such a number $(1,\pi,0,...,0)$ in hyperspherical coordinates does not invert all the vector $\vec{v}$'s components.

Also if we use the product $s \cdot \vec{v}$'s Cartesian coordinates:

$$x'' = (xx'-yy')\prod_{n=3}^{N}(1-\frac{x_n x_n'}{r_{n-1}r_{n-1}'})$$

$$y'' = (xy'+x'y)\prod_{n=3}^{N}(1-\frac{x_n x_n'}{r_{n-1}r_{n-1}'})$$

$$z'' = (r_2 z'+r_2' z)\prod_{n=4}^{N}(1-\frac{x_n x_n'}{r_{n-1}r_{n-1}'})$$

… (12.6)

$$x_k'' = (x_k r'_{k-1}+x'_k r_{k-1})\prod_{n=k+1}^{N}(1-\frac{x_n x_n'}{r_{n-1}r_{n-1}'})$$

…

$$x_N'' = x_N r'_{N-1}+x'_N r_{N-1}$$

Only $x'' = x_1''$ and $y'' = x_2''$ will be inverted.
Since the results are different this means that:
Hyperspherical $(1,\pi,0,...,0)$ and $(1,0,0,...,\pi)$ do not give the same result though they correspond to the same number in Cartesian coordinates.
One must always say: the two hypercomplex numbers are "*Cartesianly*" equivalent (equivalent in Cartesian system). Because two *Cartesianly* equivalent hypercomplex numbers may be "*hyperspherically*" different hypercomplex numbers.

This shows another intrinsic fault (flaw) inherent to the Cartesian coordinates system. Such a system does not enable us to express all geometric operations.
- It is impossible to define accurately the hypercomplex vector corresponding to scalar $(-1)$ which inverts vectors in such a coordinate system.
- It is impossible to differentiate between, for instance, different *j*'s which, when squared, give birth to different complex numbers ($j^2$ issue).

Hence we see that Cartesian coordinate system raises a lot of ambiguity when dealing with hypercomplex numbers. However, the use of hyperspherical coordinates eliminates such an equivoque, because in such a system argument definition is necessary even if the lower-order Cartesian components are null, which gives different positions corresponding to the same Cartesian position.

Example: In hyperspherical system $\vec{v}(r,\frac{\pi}{3},\frac{\pi}{2})$ is different from $\vec{v}'(r,\pi,\frac{\pi}{2})$ though they have the same Cartesian value $(0,0,r)$.

We should stop seeing scalars as dimensionless math beings (0-rank tensors, etc). Some people say they do not depend on indices or do not need any.
Indeed scalars are vectors; positive scalars have one particular property: They are colinear with *x*-axis / of the form $(r,0,...,0)$. Negative scalars are of the form $(r,0,...,\pi)$.

**Note**
$(r,0,0,...,\theta_i = \pi,0,0,...,\theta_N = 0)$ multiplies the first $i$ $x_i$-components of *N*-D vector $\vec{v}$ by $-r$.
$(r,0,0,...,\theta_i = \pi,0,0,...,\theta_k = \pi,0,0,...,\theta_N = 0)$ multiplies the $i$+1 to $k$ components of vector $\vec{v}$ by $-r$.

# 13. Applications

I expect the spherical and hyperspherical hypercomplex algebra (new vector algebra) will have an extensive and tremendous field of applications. It will have a great impact and repercussions on mathematics, physics, chemistry, biology, robotics, computer graphics and animations, etc.

In few terms, it will have a tremendous impact on all Science branches. For I strongly believe that nature itself is spherically hypercomplex.

Spherical hypercomplex fractals will be very useful in nearly all nature sciences (geology, biology, medicine, chemistry, any sort of modeling, etc.)

The multivalued $i_k^2$ may be behind quantum mechanics uncertainty principle, as well as its indeterministic and probabilistic character, because it involves a longitude and latitudes for $i_k$.

## 13.1. Spherical hypercomplex fractals: 3D Mandelbrot set
**Complex Mandelbrot Set**
The Mandelbrot set is a fractal that is defined as the set of points $c$ in the complex number plane for which the iteratively defined sequence

$$z_{n+1} = z_n^2 + c \tag{13.1}$$

with $z_0 = 0$ does not tend to infinity.[15] Fig.13.1.

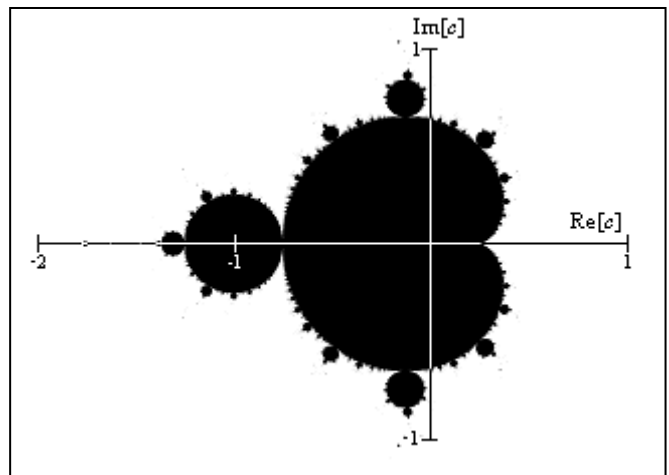

Fig.13.1. *Mandelbrot set.*[16]

In fact, we can prove that the Mandelbrot set is contained in the disk of radius 2 centered at *O*. In reality we limit the iterations by the maximum number of iterations *Nmax*. If after *Nmax* iterations $|z_{n+1}| \leq 2$ then the *c* value associated with $z_{n+1}$ is more than likely a member of the Mandelbrot set. Thus, we color *c*'s grid point black.

---

[15] Mandelbrot set, Wikipedia, the free encyclopedia, http://en.wikipedia.org/wiki/Mandelbrot_set Accessed 25 March 2014
[16] Ibidem

For $|c| \leq 2$, if we find a $z_n$ such that $|z_n| > 2$ then $c$ does not belong to the Mandelbrot set.

Using the iterative formula $h_{n+1} = h_n^2 + c$, with $h_{n+1}$, $h_n$ and $c$ hypercomplex, we can define generalized Mandelbrot sets and determine according to convergence test which numbers belong to those sets.
It will be possible to construct not only 3-D fractals, but also fractals in any *N*-D hyperspace, $N \geq 3$.

### 13.1.1. First approach

$h_0 = (0,0,0)$, $c = (x, y, z)$, $h_n = (x_n, y_n, z_n)$, $h_{n+1} = (x_{n+1}, y_{n+1}, z_{n+1})$ (13.2)

From equation (13.1) we get:
$h_{n+1} = h_n^2 + c$

$$x_{n+1} = (x_n^2 - y_n^2)\left(1 - \frac{z_n^2}{x_n^2 + y_n^2}\right) + x$$

$$y_{n+1} = 2x_n y_n \left(1 - \frac{z_n^2}{x_n^2 + y_n^2}\right) + y \quad (13.3)$$

$$z_{n+1} = 2z_n \sqrt{x_n^2 + y_n^2} + z$$

This is for $(x_n, y_n) \neq 0$,

For $n = 0$, $h_1 = c$ (13.4)

For $(x_n, y_n) = 0$ we use the product's geometric form and we need to set $\theta$ for $h_n$, we get:

$$x_{n+1} = -z_n^2 \cos 2\theta + x$$
$$y_{n+1} = -z_n^2 \sin 2\theta + y \quad (13.5)$$
$$z_{n+1} = 0 + z$$

We can set $\theta = 0$ then we get:

$$x_{n+1} = -z_n^2 + x$$
$$y_{n+1} = 0 + y = y \quad (13.6)$$
$$z_{n+1} = 0 + z = z$$

As in complex case, for each iteration we check $|h_{n+1}| = \sqrt{x_{n+1}^2 + y_{n+1}^2 + z_{n+1}^2}$. If it is greater than 2 then c does not belong to Mandelbrot set. Otherwise we continue iterations till we reach *Nmax* then we color *c*'s grid point black.

The 3D Mandelbrot fractal graphics shown in this paper have been calculated by myself early in July 2004 using DIGITAL Visual Fortran v.5.0.A and plotted on TeraPlot graphing software in May 2013 and March 2014. The maximum number of iterations is *Nmax*=100.

Figure 13.2 shows the first approach's fractal calculation results. It shows cross sections of the 3D Mandelbrot set as well as 3D plot of this set.

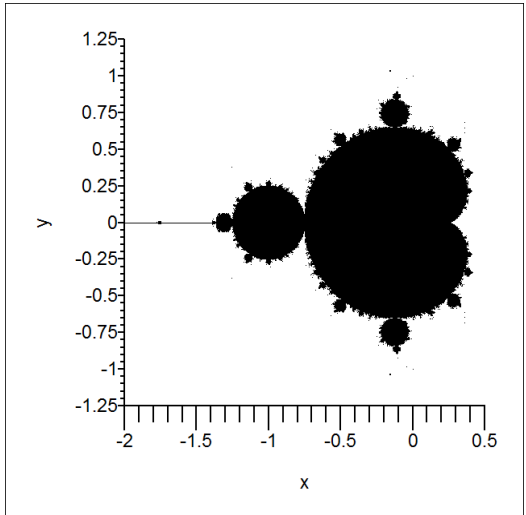
Cross section at $z = 0$

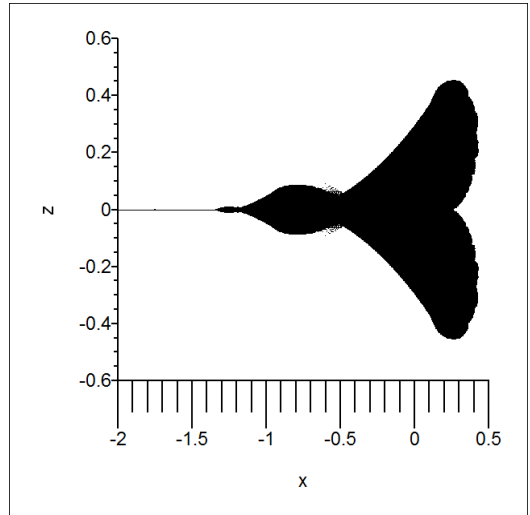
Cross section at $y = 0$

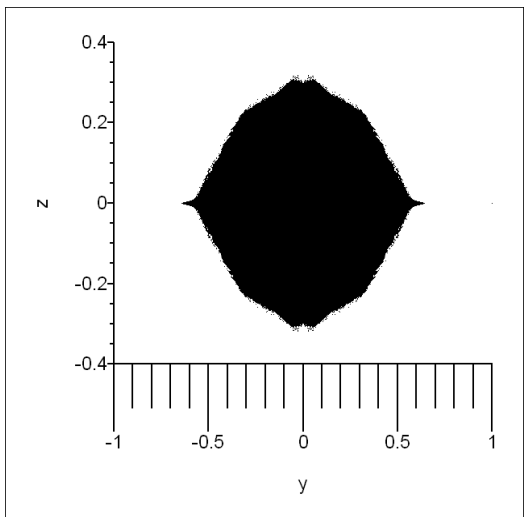
Cross section at $x = 0$

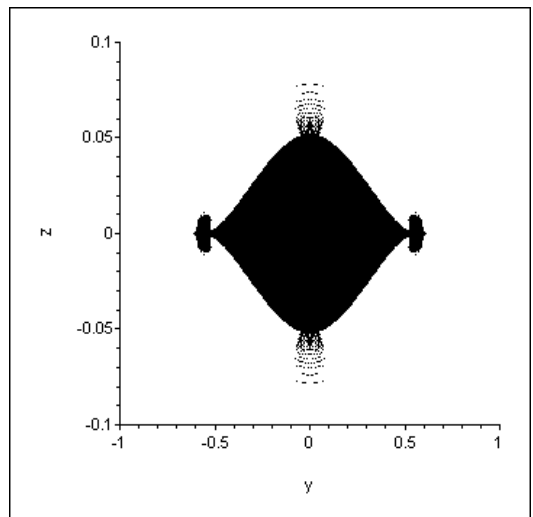
Cross section at $x = -0.5$

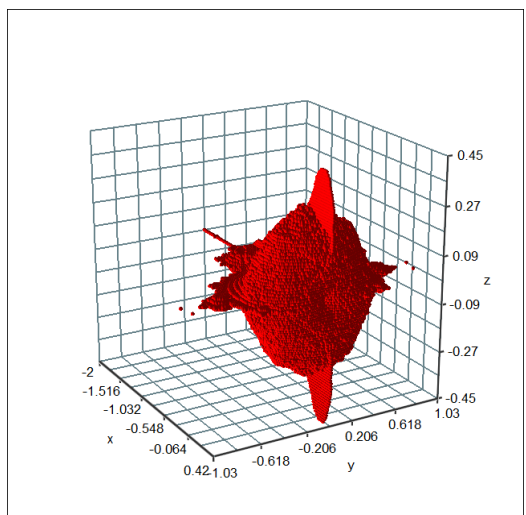
3D Mandelbrot set

Fig. 13.2. *Three-dimensional Mandelbrot set obtained with the first approach. Maximum number of iterations Nmax = 100.*

Equations (13.3) are similar to White's formula for the 3D Mandelbulb except for the $z_{n+1}$ component.[17]

### 13.1.2. Second approach

For this approach, when we calculate the Mandelbrot fractals for plane (*xOz*), that is to say *y*=0, we consider $\theta = 0$ and $\varphi$ varies from $-\pi$ to $+\pi$.

$$x = r\cos\varphi$$
$$z = r\sin\varphi \tag{13.7}$$

But in complex plane (*xOy*) we have

$$x = \rho\cos\theta$$
$$y = \rho\sin\theta \tag{13.8}$$

We notice that this plane (*xOz*) is complex and is the same as complex plane (*xOy*) with only variable change of $\theta$ into $\varphi$ and of *y* into *z* and of $\rho$ into *r*. Hence we get a conventional 2D Mandelbrot set for the plane (*xOz*).

We get the usual Mandelbrot set with the unique difference that we replace *y* with *z* and $\theta$ with $\varphi$ and $\rho$ with *r*.

$$x = r\cos\varphi$$
$$z = r\sin\varphi \tag{13.9}$$

For $y \neq 0$ when we calculate the 3D fractals we assume theta: $-\frac{\pi}{2} < \theta \leq \frac{\pi}{2}$ and $-\pi < \varphi \leq \pi$.

$$h = (x, y, z) = (r, \theta, \varphi)$$

Determination of $\theta$ and $\varphi$:

$$-\frac{\pi}{2} < \theta \leq \frac{\pi}{2}, \ -\pi < \varphi \leq \pi$$

$$\theta = \frac{\pi}{2} \text{ for } x=0$$

$$\theta = \arctan(y/x) \quad \text{for } x \neq 0$$

$$\varphi = \arctan\left(\frac{z}{\sqrt{x^2 + y^2}}\right) \text{ for } x > 0$$

$$\varphi = \pi - \arctan\left(\frac{z}{\sqrt{x^2 + y^2}}\right) \text{ for } x < 0, z \geq 0$$

$$\varphi = -\pi - \arctan\left(\frac{z}{\sqrt{x^2 + y^2}}\right) \text{ for } x < 0, z < 0 \tag{13.10}$$

$$\varphi = \arctan\left(\frac{z}{y}\right) \text{for } x = 0, y > 0$$

---

[17] Paul Nylander, Hypercomplex Fractals, http://www.bugman123.com/Hypercomplex/index.html Accessed 16 February 2014

$\varphi = \pi - \arctan\left(\dfrac{z}{|y|}\right)$ for $x = 0$ and $y < 0$ and $z \geq 0$

$\varphi = -\pi - \arctan\left(\dfrac{z}{|y|}\right)$ for $x = 0$ and $y < 0$ and $z < 0$

$\varphi = \dfrac{\pi}{2}$ for $x = 0$ and $y = 0$ and $z > 0$

$\varphi = -\dfrac{\pi}{2}$ for $x = 0$ and $y = 0$ and $z < 0$

$\varphi$ undetermined for $x = 0$ and $y = 0$ and $z = 0$

Now let's find out how to calculate the iterative formula $h_{n+1} = h_n^2 + c$.

$h_n = (x_n, y_n, z_n) = (r_n, \theta_n, \varphi_n)$ with $-\dfrac{\pi}{2} < \theta_n \leq \dfrac{\pi}{2}$ and $-\pi < \varphi_n \leq \pi$

$x_{n+1} = r_n \cos(2\theta_n)\cos(2\varphi_n) + x$
$y_{n+1} = r_n \sin(2\theta_n)\cos(2\varphi_n) + y$ \hfill (13.11)
$z_{n+1} = r_n \sin(2\varphi_n) + z$

Then we determine $h_{n+1} = (r_{n+1}, \theta_{n+1}, \varphi_{n+1})$ and so on.
Figure 13.3 shows 3D Mandelbrot set graphics for the second approach.

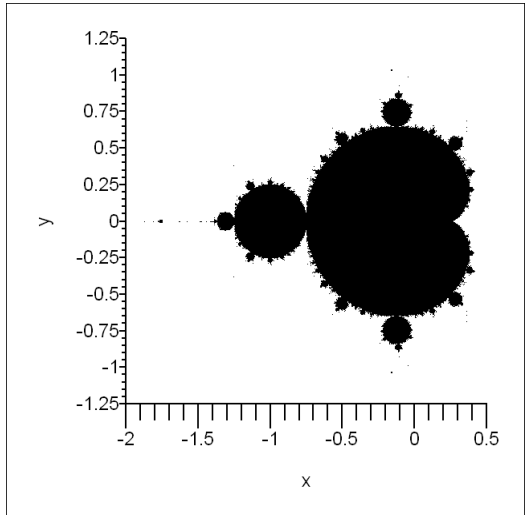
Cross section at $z = 0$

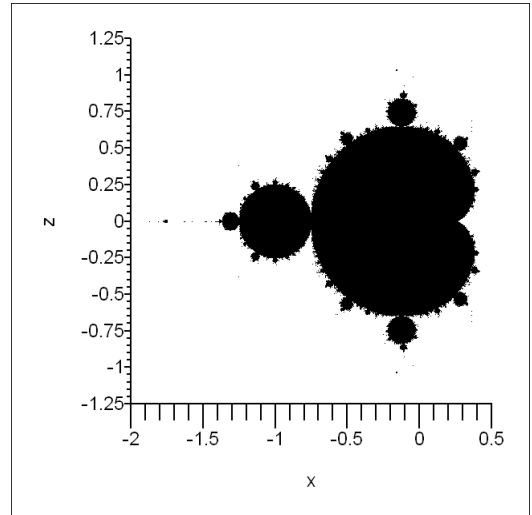
Cross section at $y = 0$

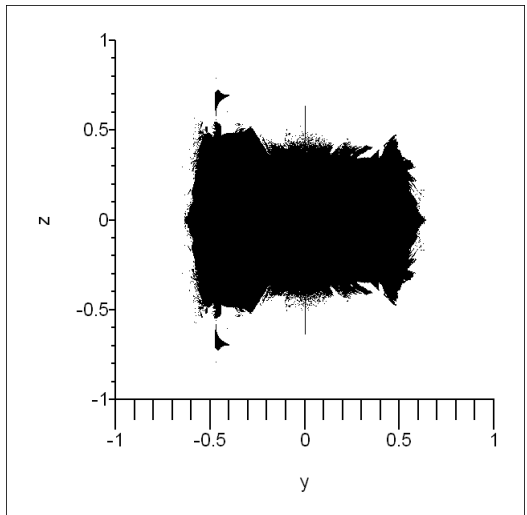
Cross section at $x = 0$

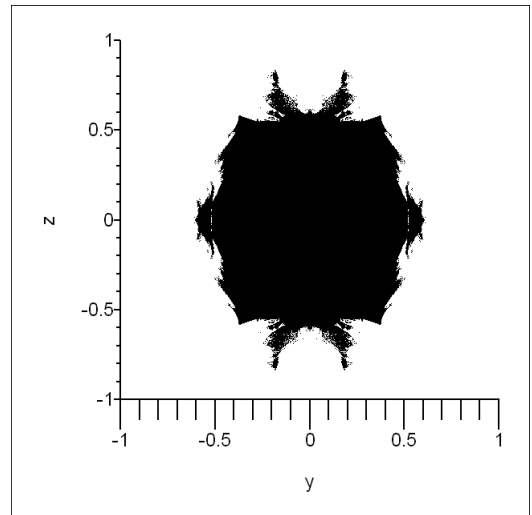
Cross section at $x = -0.5$

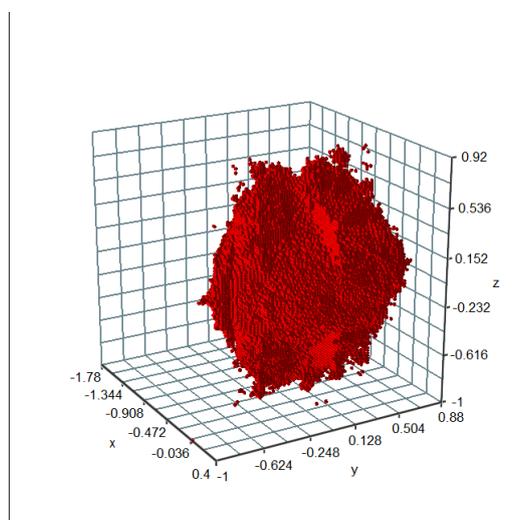
3D Mandelbrot set

Fig 13.3. *Three-dimensional Mandelbrot set obtained with the second approach. Maximum number of iterations Nmax = 100.*

## 13.2. Special and General Relativity

Hypercomplex numbers especially 4D hypercomplex numbers will certainly have broad applications in the field of special and general relativity theories, because the latter deal with 4D space-time.

**Application to the relativistic Minkowski space**

In 4D Minkowski space, let's consider two events whose coordinates are $(x_1, y_1, z_1, ct_1)$ and $(x_2, y_2, z_2, ct_2)$ in the system $S$, $c$ is the speed of light.
$ds$ is the interval between these two events

$$ds^2 = c^2(t_2 - t_1)^2 - (x_2 - x_1)^2 - (y_2 - y_1)^2 - (z_2 - z_1)^2 \qquad (13.12)$$

$S'$ is another system which is in straight motion relative to $S$ with constant speed $v$. In $S'$ the two events coordinates are $(x'_1, y'_1, z'_1, ct'_1)$ and $(x'_2, y'_2, z'_2, ct'_2)$ which are obtained through Lorentz transformation.

$$ds'^2 = c^2(t'_2 - t'_1)^2 - (x'_2 - x'_1)^2 - (y'_2 - y'_1)^2 - (z'_2 - z'_1)^2 \qquad (13.13)$$

We can prove that
$$ds^2 = ds'^2 \qquad (13.14)$$

Putting $t_2 - t_1 = dt$, $x_2 - x_1 = dx$, $y_2 - y_1 = dy$, $z_2 - z_1 = dz$ in (13.12) we get:
$ds^2 = c^2 dt^2 - (dx^2 + y^2 + dz^2)$ which is invariant to Lorentz transformation.
$ds$ is called the interval between two events or Minkowski space-time interval.[18]

Let $h = (r, \theta, \varphi, \psi) = (dx, dy, dz, dw)$ with $dw = cdt$ be the 4D hypercomplex number expressing Minkowski space-time interval.

$$r = \sqrt{dx^2 + dy^2 + dz^2 + dw^2} \qquad (13.15)$$
$$h = re^{i\theta + j\varphi + k\psi} \qquad (13.16)$$
$$h = re^{i\theta + j\varphi}(\cos\psi + k\sin\psi) \qquad (13.17)$$
$$h = re^{i\theta + j\varphi}\cos\psi + kr\sin\psi \text{ because } e^{i\theta + j\varphi}k = k \qquad (13.18)$$
$$h = r\cos\psi.e^{i\theta}(\cos\varphi + j\sin\varphi) + kr\sin\psi \qquad (13.19)$$
$$h = r\cos\psi.e^{i\theta}\cos\varphi + jr\cos\psi\sin\varphi + kr\sin\psi \qquad (13.20)$$
$$h = r\cos\psi\cos\varphi\cos\theta + ir\cos\psi\cos\varphi\sin\theta + jr\cos\psi\sin\varphi + kr\sin\psi \qquad (13.21)$$
$$h^2 = (r^2, 2\theta, 2\varphi, 2\psi) \qquad (13.22)$$
$$h^2 = r^2\cos 2\psi\cos 2\varphi\cos 2\theta + ir^2\cos 2\psi\cos 2\varphi\sin 2\theta + jr^2\cos 2\psi\sin 2\varphi + kr^2\sin 2\psi$$
$$h^2 = x + iy + jz + kw \qquad (13.23)$$

Let's find the $h^2$'s 3D component modulus.

$r_3(h^2) = \sqrt{x^2 + y^2 + z^2} = \sqrt{(r^2\cos 2\psi)^2} = |r^2\cos 2\psi| = |r^2\cos^2\psi - r^2\sin^2\psi| = |r_3^2 - dw^2|$
$= |dx^2 + dy^2 + dz^2 - c^2 t^2| = |-ds^2|$ which is invariant in Lorentz transformation.
$r_3 = \sqrt{dx^2 + dy^2 + dz^2} = r\cos\psi$ is $h$'s 3D component modulus.

**Other demonstration**

As in constructing the 3D hypercomplex numbers, the 4D hypercomplex numbers are constructed through the orthogonal complex plane construction (Fig.13.4.).
The 3D space component of $h^2$ is $r_3^2 - dw^2 = r_3^2 - c^2 dt^2 = -ds^2$

---

[18] Stamatia Mavridès (1988), Que sais-je? La relativité, Presses Universitaires de France, Paris, pp. 47-52.

This number can be negative when $\frac{\pi}{2} < \psi < \frac{3\pi}{2}$

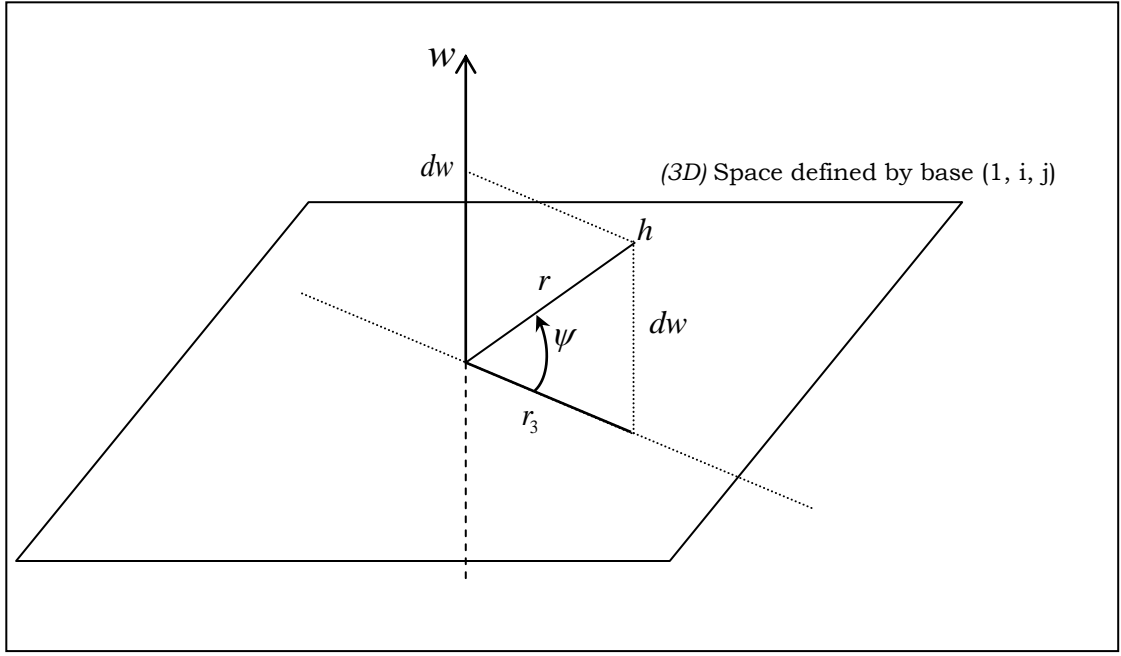

Fig. 13.4. Minkowski space and 4D hyperspherical hypercomplex numbers

$-ds^2$ = modulus of 3D spatial component of $h^2$ = 3D modulus of $h^2$ = constant for the two events.

Minkowski space doesn't take into account the *t*-axis component of the hyper number $h^2$ which is:

$$w = dh^2_{time} = 2cdt\sqrt{dx^2 + dy^2 + dz^2} \qquad (13.24)$$

$$dh^2_{time} = 2cdt\, r_3 \qquad (13.25)$$

This component may be of use in future studies of the theory of relativity.

## 14. Conclusion

The numbers described in this paper give the definitive answer to the generalization or extension of complex numbers' concept to higher dimensions. They prove to be of a wide scope of application in mathematics and science.

# 15. References

[1] Redouane Bouhannache, Solving the old problem of hypercomplex number fields: The new commutative and associative hypercomplex number algebra and the new vector field algebra; Number, Time, Relativity. Proceedings of International Scientific Meeting, Moscow: 10 – 13 August, 2004
https://sites.google.com/site/bouhannache/home/files-containing-my-contributions-on-hypercomplex-numbers Accessed  05 June 2013,
http://hypercomplex.xpsweb.com/articles/194/en/pdf/sbornik.pdf Accessed 05 June 2013
[2] Quaternion, http://www.fact-index.com/q/qu/quaternion.html Accessed 10 Jan 2014
[3] Ibid.
[4] Octonion, Wikipedia, the free encyclopedia http://en.wikipedia.org/wiki/Octonion Accessed 03 April 2004
[5] Sketching the History of Hypercomplex Numbers,
http://history.hyperjeff.net/hypercomplex Accessed 03 April 2004
[6] History of quaternions, Wikipedia, the free encyclopedia,
http://en.wikipedia.org/wiki/History_of_quaternions   Accessed 13 March 2014
[7] A Brief History of Quaternions,
http://world.std.com/~sweetser/quaternions/intro/history/history.html Accessed 03 April 2004
[8] Spherical coordinate system, Wikipedia, the free encyclopedia
http://en.wikipedia.org/wiki/Spherical_coordinate_system   Accessed 25 November 2003
[9] Field (mathematics), Wikipedia, the free encyclopedia,
http://en.wikipedia.org/wiki/Field_(mathematics) Accessed 25 November 2003
[10] Ibid.
[11] Ibid.
[12] Abelian group, Wikipedia, the free encyclopedia,
http://en.wikipedia.org/wiki/Abelian_group Accessed 25 Nov 2003
[13] Scalar, Wikipedia, the free encyclopedia
http://en.wikipedia.org/wiki/Scalar_(mathematics) Accessed 09 May 2004
[14] Ibidem
[15] Mandelbrot set, Wikipedia, the free encyclopedia,
http://en.wikipedia.org/wiki/Mandelbrot_set Accessed 25 March 2014
[16] Ibidem
[17] Paul Nylander, Hypercomplex Fractals,
http://www.bugman123.com/Hypercomplex/index.html Accessed 16 February 2014
[18] Stamatia Mavridès (1988), Que sais-je? La relativité, Presses Universitaires de France, Paris, pp. 47-52.
[19]  Eric W. Weisstein, "Hypercomplex Number" From *MathWorld*--A Wolfram Web Resource,  http://mathworld.wolfram.com/HypercomplexNumber.html Accessed 19 February 2004
[20] Abdelkader Sami et al. (1991), Mathématiques – Tome I (in Arabic), Livre scolaire de la 3$^e$ année secondaire – série mathématiques, Office National des Publications Scolaires, Algiers, pp. 113-121.

[21] S. F. Ellermeyer, The Field of Complex Numbers, pp.1-6, http://math.kennesaw.edu/~sellerme/sfehtml/classes/math4361/chapter2section6outline.pdf Accessed 21 November 2013

# CONTENTS